\documentclass[aip,
               amsmath,
               amssymb,
               floatfix,
               preprint,
               ]{revtex4-1}

\draft 

\usepackage[section]{placeins}
\def\imgwidth{0.75 \linewidth}

\usepackage[utf8]{inputenc}
\usepackage{graphicx}
\usepackage{nth}
\usepackage{physics}
\usepackage[group-separator={,}]{siunitx}
\usepackage{subcaption}
\usepackage{xcolor}

\DeclareMathOperator{\expect}{\mathbb{E}}
\DeclareMathOperator{\variance}{\mathrm{Var}}

\newcommand{\realno}[1][]{\mathbb{R}^{#1}}

\newcommand{\hp}[1]{#1_{hp}}

\def\trans{^\top} 
\def\diff{\mathrm{~d}}

\def\balpha{\boldsymbol\alpha}

\def\u{\mathbf{u}}
\def\uIC{\u_{\mathrm{IC}}}
\def\uhp{\hp{\mathbf{u}}}

\def\vs{\mathbf{v}}

\def\Jinf{J_\infty}
\def\Jref{J_\mathrm{ref}}
\def\JTh{J_{T, hp}}

\def\eTh{e_{T, hp}}

\def\dt{\Delta t}
\def\emodel{e_\mathrm{model}}

\def\Mens{M_\mathrm{ens}}

\def\eMC{e_{\mathrm{model}, \mathrm{MC}}}
\def\eMCopt{e_{\mathrm{MC}, \mathrm{opt}}}
\def\dtMC{(\dt)_{\mathrm{MC}}}
\def\dtMCopt{(\dt)_{\mathrm{MC}, \mathrm{opt}}}
\def\TsMC{T_{s, \mathrm{MC}}}

\def\NMC{N_{\mathrm{MC}}}
\def\NMCopt{N_{\mathrm{MC}, \mathrm{opt}}}

\def\attractor{\mathcal{A}}

\def\apost{a posteriori}

\def\normaldist{\mathcal{N}}
\def\halfnormaldist{\mathcal{H}}

\begin{document}

\preprint{Submitted to Physics of Fluids 22 July 2022.}

\title{Output error behavior for discretizations of ergodic, chaotic ODE systems}



\author{Cory V. Frontin}
\email[]{cfrontin@mit.edu}
\affiliation{Department of Aeronautics and Astronautics, Massachusetts Institute of Technology}

\author{David L. Darmofal}
\email[]{darmofal@mit.edu}
\affiliation{Department of Aeronautics and Astronautics, Massachusetts Institute of Technology}


\date{\today}

\begin{abstract}
The use of numerical simulation for prediction of characteristics of chaotic dynamical systems inherently involves unpredictable processes. In this work, we develop a model for the expected error in the simulation of ergodic, chaotic ODE systems, which allows for discretization and statistical effects due to unpredictability. Using this model, we then generate a framework for understanding the relationship between the sampling cost of a simulation and the expected error in the result, and explore the implications of the various parameters of simulations. Finally, we generalize the framework to consider the total cost-- including unsampled spin-up timesteps-- of simulations and consider the implications of parallel computational environments, to give a realistic model of the relationship between wall-clock time and the expected error in simulation of a chaotic ODE system.
\end{abstract}


\maketitle 



\section{Introduction}

For chaotic systems, estimation of long-time behavior is challenging because chaotic systems have limited predictability \cite{lighthill1986recently}. Of the general class of chaotic systems, a subset are ergodic systems, whose long-term states are drawn from a stationary distribution, independent of initial condition \cite{eckmann1985ergodic}. For ergodic chaotic problems, we frequently want to quantify the unique infinite-time average of some instantaneous quantity of interest of the system:
\begin{equation}
    \Jinf= \lim_{T \to \infty} \frac{1}{T} \int_{0}^{T} g(\u(t)) \diff t ,
    \label{eq:output_true}
\end{equation}
where $g$ is the instantaneous output functional, and the state $\u(t)$ is governed by a dynamical system of the form:
\begin{equation}
    \dv{\u}{t}= f(\u)
\end{equation}
with a given initial condition (IC), $\u(0)= \uIC$.

Often, the complexity of a chaotic systems of interest is high, and accordingly the cost of an accurate computational estimate of $\Jinf$ becomes formidable \cite{chapman1979computational, spalart1997comments, choi2012grid}. As the cost of computational simulation gets larger, efficient discretization methods become critical for accurately estimating quantities of interest.

Understanding the error in approximations of $\Jinf$ is nontrivial because statistical errors (errors due to finite-time approximation) and discretization error (error due to numerical approximation of solutions) are always simultaneously present. In the largest Direct Numerical Simulation (DNS) and Large Eddy Simulation (LES) cases, for example, it is typical to fix sampling time at some large number of characteristic times and validate that discretization error converges as expected, assuming negligible sampling error \cite{kim1987turbulence, lozano2014effect, delalamo2004scaling, goc2021large}. Recent work has sought to quantify the effect of statistical error more robustly, using turbulent flow theory \cite{thompson2016methodology}, advanced spatio-temporal statistical post-processing methods \cite{russo2017fast}, statistical windowing techniques \cite{mockett2010detection}, or by extending the concept of Richardson extrapolation to chaotic flows using auto-regressive models and Bayesian methods \cite{oliver2014estimating}. The latter work is notable for its use to estimate the statistical errors in the DNS of a high-$\Re$ turbulent channel flow \cite{lee2015direct}.

The objective of this paper is to investigate the behavior of statistical and discretization errors as a function of computational cost for ergodic systems. Following a similar approach to \citet{oliver2014estimating}, we propose a simple error model for finite-time, discrete approximations of infinite-time averages on attractors. Using the Lorenz system as an example, we demonstrate that the discretization error converges as timestep size decreases. However, it does not increase exponentially with sampling time as might be expected from classical numerical analysis but rather asymptotes to a constant value with respect to sampling time.  Further, for a given computational cost (e.g. number of timsteps), an optimal choice of discretization (i.e. timestep) exists that minimizes the expected error in a simulation, when accounting for both the effects of discretization error and sampling error. We show that this optimal choice results in a convergence rate with respect to computational cost that is bounded by the sampling convergence rate with a minor impact from the discretization order of accuracy. Finally, we consider the implications of spin-up time (i.e. unsampled time needed to arrive at the stationary distribution) and parallelism on the optimal error. We develop a method for estimating transient-related errors, and then evaluate optimal choices incorporating the results.


\section{Proposed error model on the attractor}

To approximate $\Jinf$, we compute finite-time, discrete estimates of the outputs of interest of the true system:
\begin{equation}
    \JTh= \frac{1}{T_s} \hbox{\Large I}_{t_0}^{t_0 + T_s} \big( \hp{g}(\uhp(t)) \big),
    \label{eq:output_discrete}
\end{equation}
where the notation $\mathrm{I}_{a}^{b} (\cdot)$ here represents the quadrature approximation of the integral $\int_{a}^{b} (\cdot) \diff t$ of a quantity $(\cdot)$ between $a$ and $b$. Here, we have made a discrete approximation of the state using an order-$p$ discretization with a temporal grid with characteristic size $h= \Delta t$, where an order-$p$ discretization is one for which the discretization error behaves as:
\begin{equation}
  \max_{t \in [0, t_0 + T_s]} \abs{\hp{g}(\uhp(t)) - g(\u(t))}= \mathcal{O}(h^p)
\end{equation}
when the discretization is applied to a well-posed (non-chaotic) system. Then we sample that discrete state over a finite sampling period, $T_s$, starting at some initial time $t_0$. We can define the error that is incurred as
\begin{equation}
    e_{T, hp}= \JTh - \Jinf .
    \label{eq:error_total_def}
\end{equation}
By introducing a third value,
\begin{equation}
    J_T= \frac{1}{T_s} \int_{t_0}^{t_0 + T_s} g(\u(t)) \diff t ,
    \label{eq:output_ficticious}
\end{equation}
we can re-write the error using an identity:
\begin{equation}
    e_{T, hp}= (\JTh - J_T) + (J_T - \Jinf)= \hp{e} + e_T .
    \label{eq:error_total}
\end{equation}
Here, we define the ``discretization error'' and ``sampling error'', respectively:
\begin{align}
    \hp{e} & \equiv \JTh - J_T
    \label{eq:error_disc} \\
    e_T & \equiv J_T - \Jinf .
    \label{eq:error_stat}
\end{align}

We can take an absolute value of both sides of \eqref{eq:error_total}, followed by a manipulation using the triangle inequality:
\begin{equation}
    \begin{aligned}
        |e_{T, hp}| & = |\hp{e} + e_T| \\
        & \leq |\hp{e}| + |e_T| .
    \end{aligned}
    \label{eq:error_total_triangle}
\end{equation}
Thus, the total error incurred by approximation is bounded by the sum of the absolute discretization and sampling errors. Next, we define the attractor of the operator $f$, $\attractor$, as the set of long-term states towards which all trajectories converge independently of initial condition \cite{stuart1994numerical}. We can define the expectation $\expect_{\attractor}[\phi(\u_0)]$ for a generic function $\phi$ as the expectation taken over all the trajectories that can result from starting from points on the attractor, $\attractor$:
\begin{equation}
    \expect_{\attractor}[\phi]= \frac{1}{|\attractor|} \int_{\u_0 \in \attractor} \phi(\u_0) \diff \u_0 .
    \label{eq:expectation_attractor}
\end{equation}
For the case in question we will be considering either
\begin{equation*}
    \phi(\u_0)= \frac{1}{T_s} \left| \int_{t_0}^{t_0 + T_s} g(\u(t)) \diff t - \Jinf \right| ,
\end{equation*}
or
\begin{equation*}
    \phi(\u_0)= \frac{1}{T_s} \left| \hbox{\Large I}_{t_0}^{t_0 + T_s} \big( \hp{g}(\uhp(t)) \big) - \int_{t_0}^{t_0 + T_s} g(\u(t)) \diff t \right| ,
\end{equation*}
with, for these examples, $\u(t_0)= \u_0 \in \attractor$. Given these definitions, we can now take the expectation of \eqref{eq:error_total_triangle}, giving
\begin{align}
    \expect_{\attractor}[|e_{T, hp}|] & \leq \expect_{\attractor}[|\hp{e}|] + \expect_{\attractor}[|e_T|]
    \label{eq:error_total_expectation}
\end{align}
by linearity.

From here, we propose asymptotic forms for the two right-hand side terms in \eqref{eq:error_total_expectation}. Consider the definition of $e_T$ in \eqref{eq:error_stat}:
\begin{equation}
  e_T= \frac{1}{T_s} \int_{t_0}^{t_0 + T_s} g(\u(t)) \diff t - \Jinf .
\end{equation}
Assuming that we choose $t_0$ such that each $\u_0$ is effectively an independent sample from the attractor's stationary distribution, then the quantity $g(\u(t))$ is a random variable drawn from a stationary distribution. The states of ergodic systems, in general, are not independent in time, but as long as the system has satisfactorily strong mixing properties, the central limit theorem (CLT) can be applied to finite time averages of its outputs. This is the case whenever the condition of $\alpha$-mixing is met \cite{denker1989central, bradley2005basic}, which has been shown for the Lorenz system \cite{araujo2015rapid}. Thus we can write $e_T$ as:
\begin{equation}
  e_T \sim \mathcal{N} \left( 0, \left( \sqrt{\frac{\pi}{2}} A_0 T_s^{-1/2} \right)^2 \right) ,
\end{equation}
where $\normaldist(\mu, \sigma^2)$ gives the normal distribution with mean $\mu$ and variance $\sigma^2$. If we take the absolute value of this random variable, the result is a halfnormal distribution:
\begin{equation}
  \left| e_T \right| \sim \halfnormaldist \left( \left( \sqrt{\frac{\pi}{2}} A_0 T_s^{-1/2} \right)^2 \right) ,
  \label{eq:abs_stat_error_halfnormal}
\end{equation}
where $\halfnormaldist(\sigma^2)$ gives a halfnormal distribution such that $|X| \sim \halfnormaldist(\sigma^2)$ when $X \sim \normaldist(0, \sigma^2)$. The expectation of the half-normal distribution is well defined, allowing:
\begin{equation}
    \expect_{\attractor}[|e_T|] \approx A_0 T_s^{-1/2}
    \label{eq:statistical_error}
\end{equation}
as $T_s$ goes to infinity.

Now consider the use of a time-stepping method to give a discrete approximation $\uhp(t_n)$ of $\u(t_n)$ for each $t_n= n (\dt)$. Following classical analysis \cite{hairer1993solving}, we might expect that the discretization error should take a form:
\begin{equation}
  |\hp{e}| \approx C_p \left( \frac{\exp(\Lambda T_s) - 1}{\Lambda} \right) (\dt)^p .
\end{equation}
This analysis is based on bounding the growth of local truncation error at each timestep by the Lipschitz constant, $\Lambda$, of the underlying system, with $C_p$ a constant parameter that depends on the choice of method. However, Viswanath showed \cite{viswanath2001global} that, the global error could be modeled by a form:
\begin{equation}
  |\hp{e}| \approx E(T_s; p) (\dt)^p ,
\end{equation}
where $E(T_s; p)$ could be bounded by a constant for some nonlinear but non-chaotic systems that are exponentially stable. While this result has not been extended to an ergodic system, the expected convergence onto the attracting set suggests a bound of the form:
\begin{equation}
    \expect_{\attractor}[|\hp{e}|] \approx C_p (\dt)^{p} ,
    \label{eq:truncation_error}
\end{equation}
As our results in Section~\ref{sec:evaluation} will show, \eqref{eq:truncation_error} is a good description of the expected discretization error.

Thus, taking \eqref{eq:error_total_expectation}, \eqref{eq:statistical_error}, and \eqref{eq:truncation_error} we assume a bound of the form:
\begin{equation}
    \expect_{\attractor}[|e_{T, hp}|] \leq \emodel= C_q (\dt)^q + A_0 T_s^{-r} ,
    \label{eq:error_model}
\end{equation}
that bounds $\expect_{\attractor}[|e_{T, hp}|]$ when $\dt$ is small enough and $T_s$ is large enough to satisfy the asymptotic assumptions. Here, $q$ is the observed discretization convergence rate, which in practice may differ from $p$ due to numerical cancellations or if the solutions of the system are insufficiently regular. Similarly, $r$ is an observed sampling convergence rate coefficient, which we expect to be $1/2$ asymptotically under the CLT.


\section{Evaluation of proposed error model on the Lorenz system}
\label{sec:evaluation}

In the following section, we will fit numerical results for the Lorenz system to determine $q$, $r$, $C_q$, and $A_0$ and show that this model is representative of the observed behavior. The Lorenz system is given by \cite{lorenz1963deterministic}:
\begin{equation}
    \dv{\u}{t}= f(\u; \balpha)=
    \begin{pmatrix}
        \alpha_0 (u_1 - u_0) \\
        u_0 (\alpha_1 - u_2) - u_1 \\
        u_0 u_1 - \alpha_2 u_2
    \end{pmatrix} ,
\end{equation}
where $\u= [u_0, u_1, u_2] \trans$ and $\balpha= [\alpha_0, \alpha_1, \alpha_2] \trans$. The Lorenz system is known to be chaotic for the classic Lorenz parametrization \cite{sparrow1982lorenz}: $\balpha= [10, 28, 8/3]$, which is used everywhere in this text. For the output, we choose $g(\u)= u_2$. We consider a set of explicit methods: forward Euler (FE, $p= 1$), \nth{3}-order Runge-Kutta (RK3, $p= 3$), and \nth{4}-order Runge-Kutta (RK4, $p= 4$). In all of these methods, we expect asymptotic convergence of $\JTh$ to $J_T$ to be at least $\mathcal{O}(\dt^p)$ for non-chaotic systems \cite{dormand1984global}.

For any given discrete instance, we will start the simulation at an initial state at $t= 0$ that is sampled randomly from a normal distribution:
\begin{equation}
    \begin{aligned}
        \u_{\mathrm{init}} & \sim
        \begin{pmatrix}
            \mathcal{N}(1.0, 5.0^2) \\
            \mathcal{N}(1.0, 5.0^2) \\
            \mathcal{N}(1.0, 5.0^2)
        \end{pmatrix}
    \end{aligned} .
    \label{eq:random_ic}
\end{equation}
To guarantee that the initial sampling state $\u_0$ at $t_0$ is on the attractor (as well as further guaranteeing the independence from the other Monte Carlo instances), we evolve the state of any given Lorenz system discretization from its starting state $\u_{\mathrm{init}}$ for $t_0= 100$ before proceeding to sample; we refer to the process of evolving the solution until it is on the attractor as ``spin-up''. Then, we evolve the state over the next $T_s$, during which we integrate and compute \eqref{eq:output_discrete} using the same numerical integration scheme that was used for the state itself.

To approximate $\eTh$, we must first estimate $\Jinf$ by a reference value $\Jref$. $\Jref$ is calculated using an ensemble mean of $\JTh$ over $\Mens= 512^2$ instances of the Lorenz system. Each instance is started from a different $\u_{\mathrm{init}}$ as given in \eqref{eq:random_ic} and simulated using RK4 with $\dt= \num{17.7e-6}$ and $T_s= 6646.9$. The resulting $\Jref$ is:
\begin{equation}
    \Jref= 23.549916 \pm 0.000074 ,
    \label{eq:output_reference}
\end{equation}
with a 95\% confidence estimate based on the ensemble mean estimator.

The computation of $\Jref$ allows us to estimate errors $e_{T, hp} \approx \JTh - \Jref$. For a given $\dt$, $T_s$ pair, we then approximate $\expect_{\attractor}[|e_{T, hp}|]$ using a Monte Carlo method over $M= \num{10000}$ independent instances of the discrete system, each started from initial states drawn from \eqref{eq:random_ic} and spun-up to independent sampling starting points on the attractor $\u_0^{(m)}$:
\begin{equation}
    \expect[|e_{T, hp}|] \approx \frac{1}{M} \sum_{m= 1}^{M} \left| \JTh \left( \u_{0}^{(m)} \right) - \Jref \right| .
\end{equation}

\begin{figure}[h]
    \centering%
    \includegraphics[width= \imgwidth]{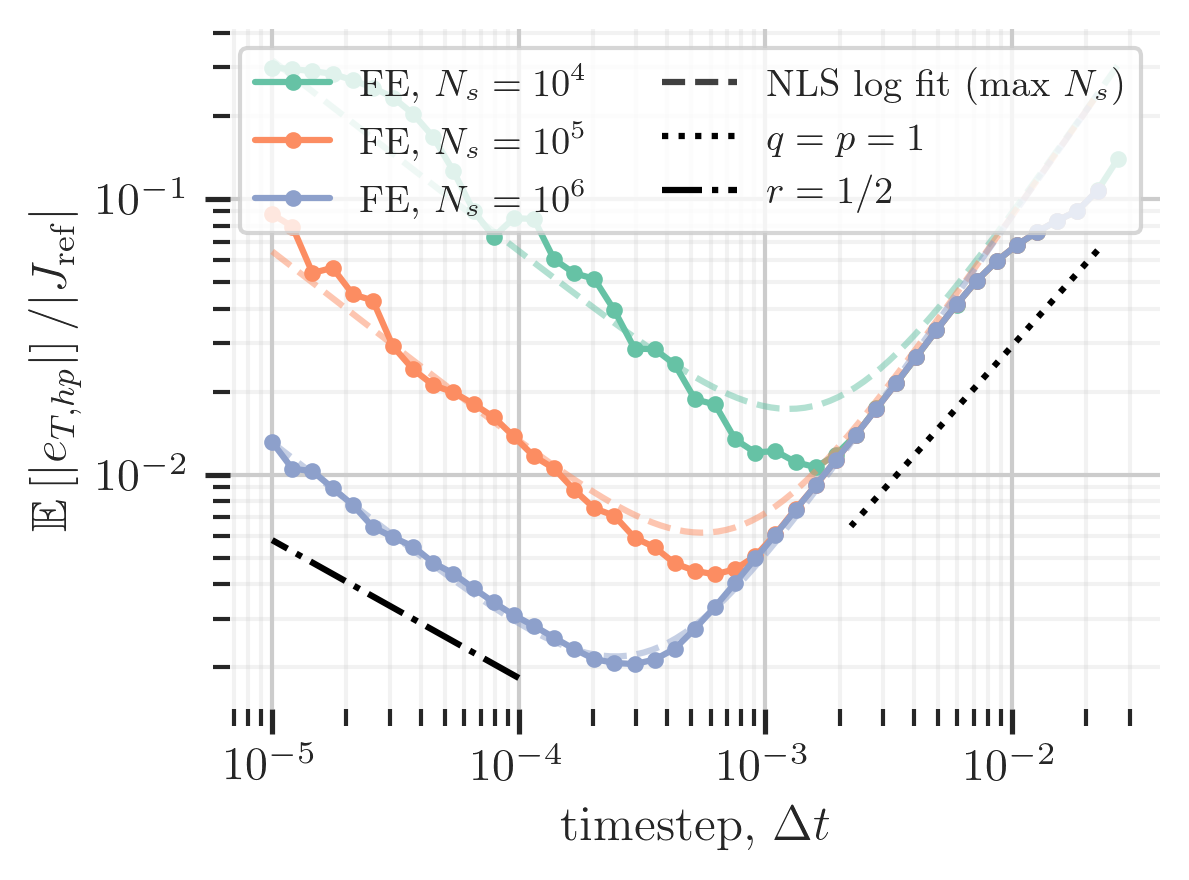}%
    \caption{Expected relative error as a function of $\dt$ for Forward Euler discretization of the Lorenz equations. Nonlinear least squares fit based on $N_s= 10^6$ data.}%
    \label{fig:errormodel_trade_FE_N1000000}%
\end{figure}

\begin{figure}[h]
    \centering%
    \includegraphics[width= \imgwidth]{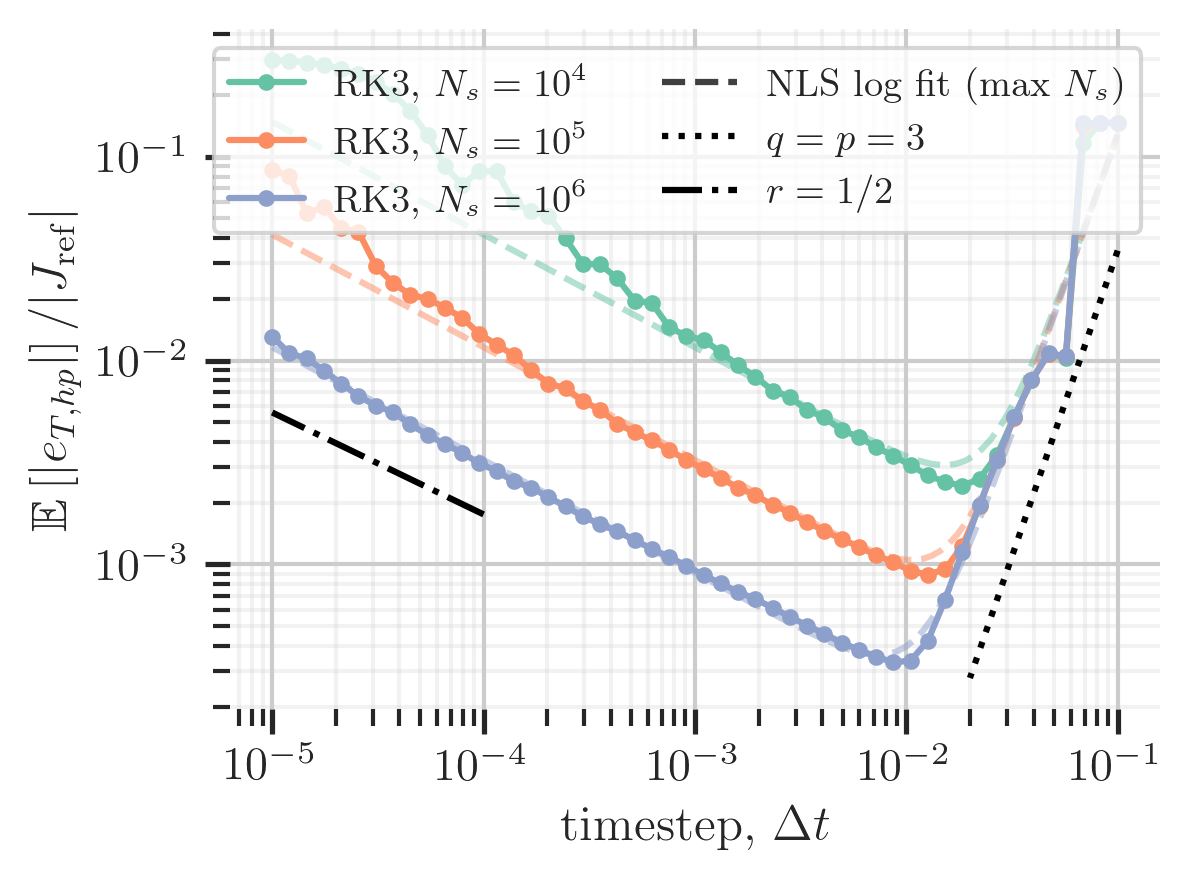}%
    \caption{Expected relative error as a function of $\dt$ for RK3 discretization of the Lorenz equations. Nonlinear least squares fit based on $N_s= 10^6$ data.}%
    \label{fig:errormodel_trade_RK3_N1000000}%
\end{figure}

\begin{figure}[h]
    \centering%
    \includegraphics[width= \imgwidth]{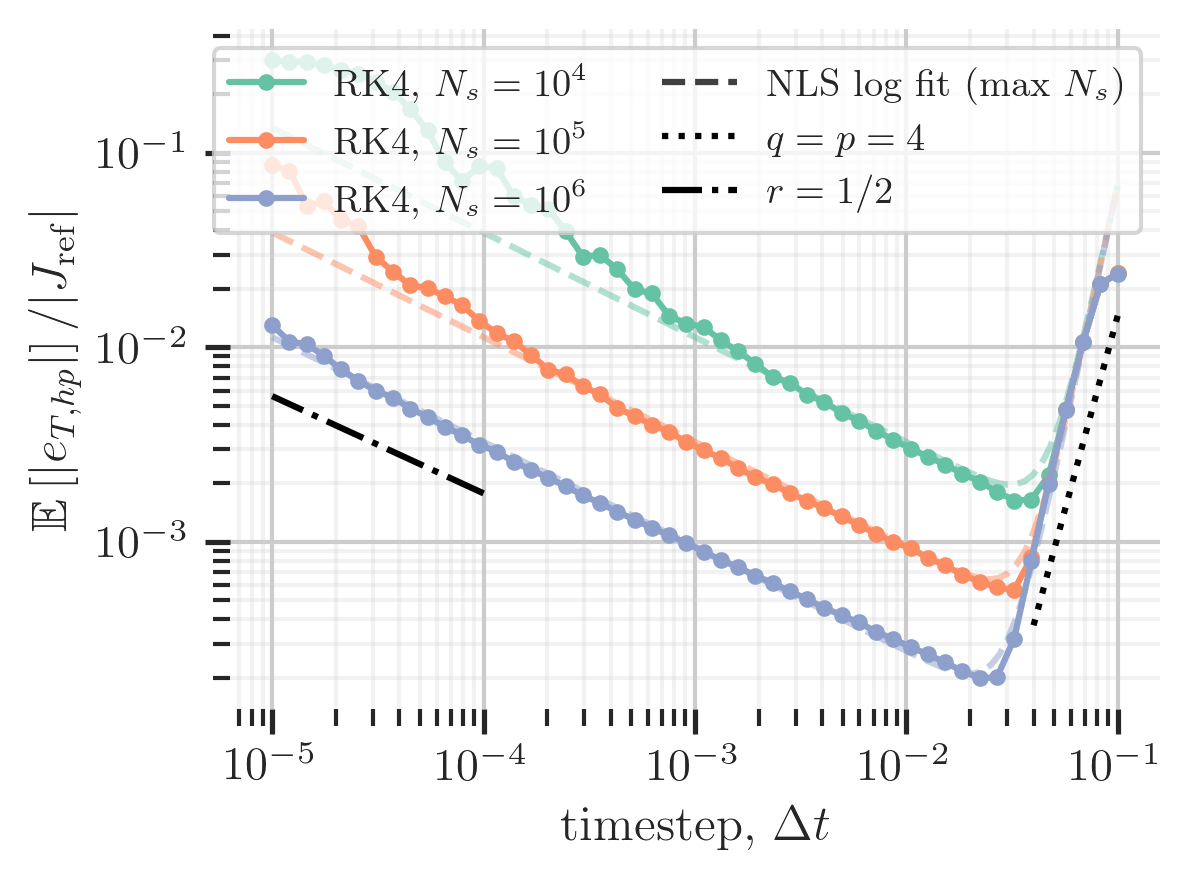}%
    \caption{Expected relative error as a function of $\dt$ for RK4 discretization of the Lorenz equations. Nonlinear least squares fit based on $N_s= 10^6$ data.}%
    \label{fig:errormodel_trade_RK4_N1000000}%
\end{figure}

In Figures~\ref{fig:errormodel_trade_FE_N1000000}, \ref{fig:errormodel_trade_RK3_N1000000},~and~\ref{fig:errormodel_trade_RK4_N1000000}, we compare the results of simulations with the FE, RK3, and RK4 discretizations with different values of $N_s$. In these figures, $T_s$ scales with $\dt$ for a given $N_s$, so the $T_s$ values on the x-axis will vary between lines on the plot. The fits shown are computed with truncated data, in order to eliminate non-convergent data at small $T_s$ or large $\dt$; the limits used for truncation are found in Table~\ref{tab:fit_bounds}. The results of the nonlinear least squares fits for $N_s= 10^4$, $10^5$, and $10^6$, are given in Table~\ref{tab:NLS_fits}. In the table, we observe that $r \to 1/2$ as the discretization error is reduced, either by increasing $N_s$ or by pushing $p$ higher.
\begin{table}
    \centering
    \begin{tabular}{l | r | r}
        method & $\dt_\mathrm{max}$ & $T_{s, \mathrm{min}}$ \\
        \hline
        FE & $5.0 \times 10^{-3}$ & $1.0$ \\
        RK3 & $5.0 \times 10^{-2}$ & $1.0$ \\
        RK4 & $9.0 \times 10^{-2}$ & $1.0$
    \end{tabular}
    \caption{Fit boundaries for nonlinear least squares fits.}
    \label{tab:fit_bounds}
\end{table}
These figures demonstrate that \eqref{eq:truncation_error} has explanatory value, as the errors in the discretization-dominated region collapse independently of $T_s$. It is also worth noting that Table~\ref{tab:NLS_fits} demonstrates higher-than-expected discretization error convergence rates for FE and RK4.
\begin{table}
    \centering
    \begin{subtable}{0.495 \linewidth}
      \centering
      \begin{tabular}{l | r | r | r}
          & FE & RK3 & RK4 \\
          \hline
          $A_0$ & {2.19} &    {1.74} &    {1.63} \\
          $r$ &   {0.975} &   {0.721} &   {0.683} \\
          $C_q$ & {4995} &    {942} &     \num{85900} \\
          $q$ &   {1.65} &    {2.70} &    {4.83} \\
      \end{tabular}
      \subcaption{$N_s= 10^4$}
      \label{tab:NLS_fits_tiny}
    \end{subtable}
    \begin{subtable}{0.495 \linewidth}
      \centering
      \begin{tabular}{l | r | r | r}
          & FE & RK3 & RK4 \\
          \hline
          $A_0$ & {1.94} &    {1.50} &    {1.41} \\
          $r$ &   {0.820} &   {0.648} &   {0.620} \\
          $C_q$ & {1410} &    {1310} &    \num{96100} \\
          $q$ &   {1.40} &    {2.76} &    {4.84} \\
      \end{tabular}
      \subcaption{$N_s= 10^5$}
      \label{tab:NLS_fits_small}
    \end{subtable}%
    \hfill%
    \begin{subtable}{0.495 \linewidth}
        \centering
        \begin{tabular}{l | r | r | r}
            & FE & RK3 & RK4 \\
            \hline
            $A_0$ & {1.52} &    {0.978} &   {0.918} \\
            $r$ &   {0.693} &   {0.553} &   {0.538} \\
            $C_q$ & {714.6} &   {2740} &    \num{165000} \\
            $q$ &   {1.273} &   {2.96} &    {5.02} \\
        \end{tabular}
      \subcaption{$N_s= 10^6$}
      \label{tab:NLS_fits_big}
    \end{subtable}
    \caption{Values of error model coefficients computed from nonlinear least squares fits to Monte Carlo study data.}
    \label{tab:NLS_fits}
\end{table}%
\begin{figure}[h]
    \centering%
    \includegraphics[width= \imgwidth]{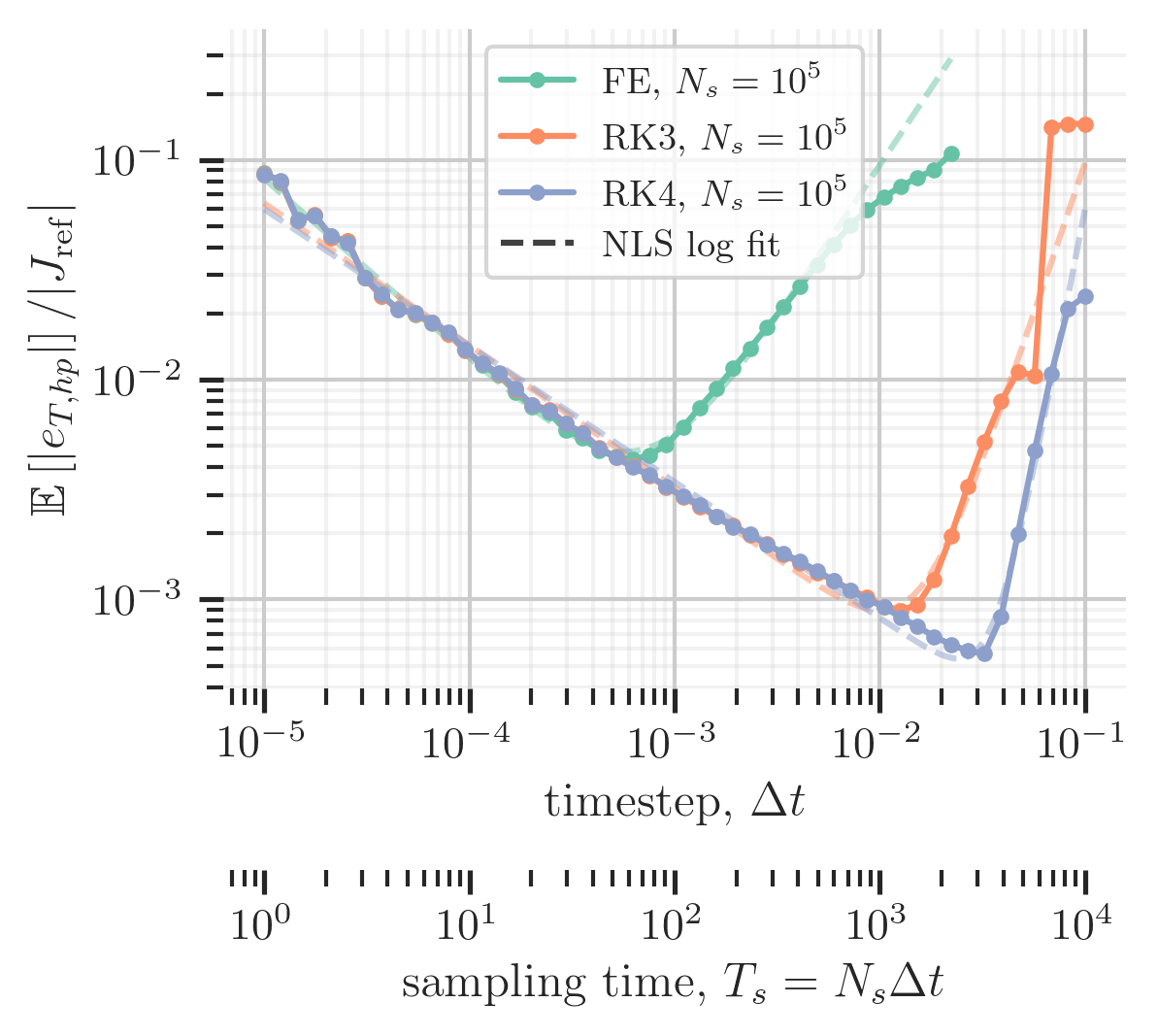}%
    \caption{Expected relative error as a function of $\dt$ for discretizations of the Lorenz equations.}%
    \label{fig:errormodel_trade_N100000}%
\end{figure}
In Figure~\ref{fig:errormodel_trade_N100000}, we can examine the sampling error behavior between discretization methods for a single shared choice of $N_s$. Here, we can see that the sampling error effects on the left-hand side of the plot collapse independently of the discretization method. This indicates that the statistical effects are properties of the dynamical system, not artifacts of the discretization, as we might expect in the limit as $\dt \to 0$.

Finally, we attempt to compare the computational costs across the various discretizations. In this case, the number of timesteps $N_s$ is not a good proxy for fixed cost, since the computation time for a timestep will vary between methods. Instead, we now fix $U_s$, the total number of evaluations of the right-hand side $f$ used in sampling timesteps. For the explicit schemes used in this work, we will have $p$ right-hand side evaluations (e.g. Forward Euler has $p= 1$ right-hand side evaluations), and thus $U_s= p N_s$.
\begin{figure}
    \centering
    \includegraphics[width= \imgwidth]{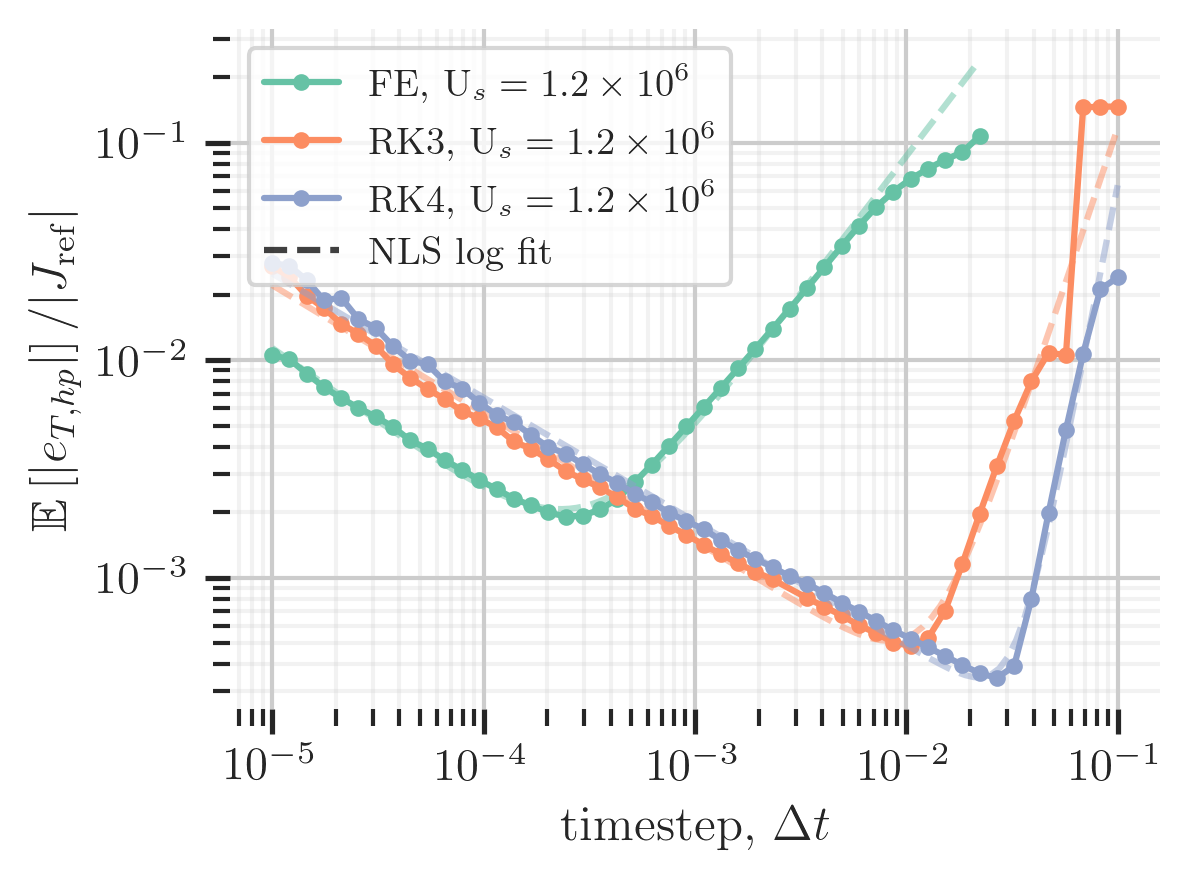}
    \caption{Expected percent error as a function of $\dt$ for discretizations of the Lorenz equations at a number of sampling residual evaluations. All fits evaluated at $U_s= \num{1.2e6}$.}
    \label{fig:errormodel_trade_evals}
\end{figure}
In Figure~\ref{fig:errormodel_trade_evals}, we can see the effect of changing $\dt$ at fixed sampling cost $U_s$ across discretizations. The error that can be achieved with the Runge-Kutta methods is lower than that of the forward Euler scheme, a factor of $4.8$ improvement in the error from FE to RK4. However, the best-case improvement for going from \nth{3}-order to \nth{4}-order Runge-Kutta schemes is a only factor of about $1.4$. Moreover, the results show that to achieve the lowest possible error, the optimal timestep will be discretization dependent. We investigate this further in the next section.


\section{Optimal timestepping on the attractor}

We now study the implications of the error model \eqref{eq:error_model}, specifically seeking to understand the convergence of the error with respect to computational effort.  In this analysis, we will assume that $r= 1/2$.

Consider a non-dimensional form of error model in which the error is normalized by the standard deviation of the instantaneous output $\sigma_g$ and the timescales $\Delta t$ and $T_s$ are normalized by decorrelation time $T_d$. The decorrelation time relates the amount of variance from independent draws from the distribution on the attractor and the amount of variance in the finite-time mean estimators based on the correlated output signal, given by the relation \cite{trenberth1984some}:
\begin{equation}
  \variance[J_T]= \frac{T_d}{T_s} \sigma_g^2 .
  \label{eq:decorr_variance_relation}
\end{equation}
Furthermore, combining \eqref{eq:abs_stat_error_halfnormal} and \eqref{eq:decorr_variance_relation} allows us to write
\begin{equation}
  A_0= \sqrt{\frac{2}{\pi}} \sigma_g T_d^{1/2} .
  \label{eq:coefficient_statistical_relation}
\end{equation}
In general, $T_d$ is hard to estimate accurately; this is a crux of the work of \citet{oliver2014estimating}. In our formulation of the error model, we identify $A_0$, which avoids outright estimation of $T_d$. However, for the purposes of understanding the behavior of the error, $T_d$ is an intrinsic timescale which can be used to normalize $\dt$ and $T_s$.

The resulting non-dimensional form of the error model is
\begin{equation}
  \frac{\emodel}{\sigma_g}= \frac{C_q T_d^q}{\sigma_g} \left( \frac{\dt}{T_d} \right)^q + \sqrt{\frac{2}{\pi}} \left( \frac{T_s}{T_d} \right)^{-\frac{1}{2}} .
  \label{eq:error_model_nondim}
\end{equation}
\begin{widetext}    
We can also write the optimizers and optimal value of \eqref{eq:error_model_nondim} in terms of the non-dimensional variables. These are given by:
\begin{equation}
  \begin{aligned}
    \left( \frac{\dt}{T_d} \right)_{\mathrm{opt}} &= \left( \frac{1}{2 \pi} \right)^{\frac{1}{2q+1}}
            \left( \frac{q C_q T_d^q}{\sigma_g} \right)^{-\frac{2}{2q+1}} N_s^{-\frac{1}{2q+1}} \\
    \left( \frac{T_s}{T_d} \right)_{\mathrm{opt}} &= \left( \frac{1}{2 \pi} \right)^{\frac{1}{2q+1}}
            \left( \frac{q C_q T_d^q}{\sigma_g} \right)^{-\frac{2}{2q+1}} N_s^{\frac{2q}{2q+1}} \\
    \left( \frac{\emodel}{\sigma_g} \right)_{\mathrm{opt}} &=
            \left( \frac{1}{2 \pi} \right)^{\frac{q}{2q+1}} \left( 2 + \frac{1}{q} \right) \left( \frac{q C_q T_d^q}{\sigma_g} \right)^{\frac{1}{2q+1}} N_s^{-\frac{q}{2q+1}} .
  \end{aligned}
  \label{eq:error_optimal_nondim}
\end{equation}
\end{widetext}

In terms of convergence with respect to sampling costs, the error model will scale at best as
\begin{equation*}
  \left( \frac{\emodel}{\sigma_g} \right)_{\mathrm{opt}} \sim N_s^{-\frac{q}{2q+1}} .
\end{equation*}
In the limit as $q \to \infty$, the rate $q/(2q+1) \to 1/2$: the CLT limits the convergence rate. Table~\ref{tab:sampling_convergence_rates} gives the rates of convergence \eqref{eq:error_optimal_nondim} for various values of $q$.

\begin{table}[h]
  \centering
  \begin{tabular}{c|c|c|c|c|c|c|c}
    $q$ &               $1$ &     $2$ &     $3$ &     $4$ &     $5$ &     $\cdots$ &  $\infty$ \\
    \hline
    $\frac{q}{2q+1}$ &  $1/3$ &   $2/5$ &   $3/7$ &   $4/9$ &   $5/11$ &  $\cdots$ &  $1/2$
  \end{tabular}
  \caption{Convergence rates for combined error with respect to sampling timesteps implied by \eqref{eq:error_optimal_nondim} at common high-order discretization error convergence rates.}
  \label{tab:sampling_convergence_rates}
\end{table}

\begin{figure}[h]
    \centering%
    \includegraphics[width= \imgwidth]{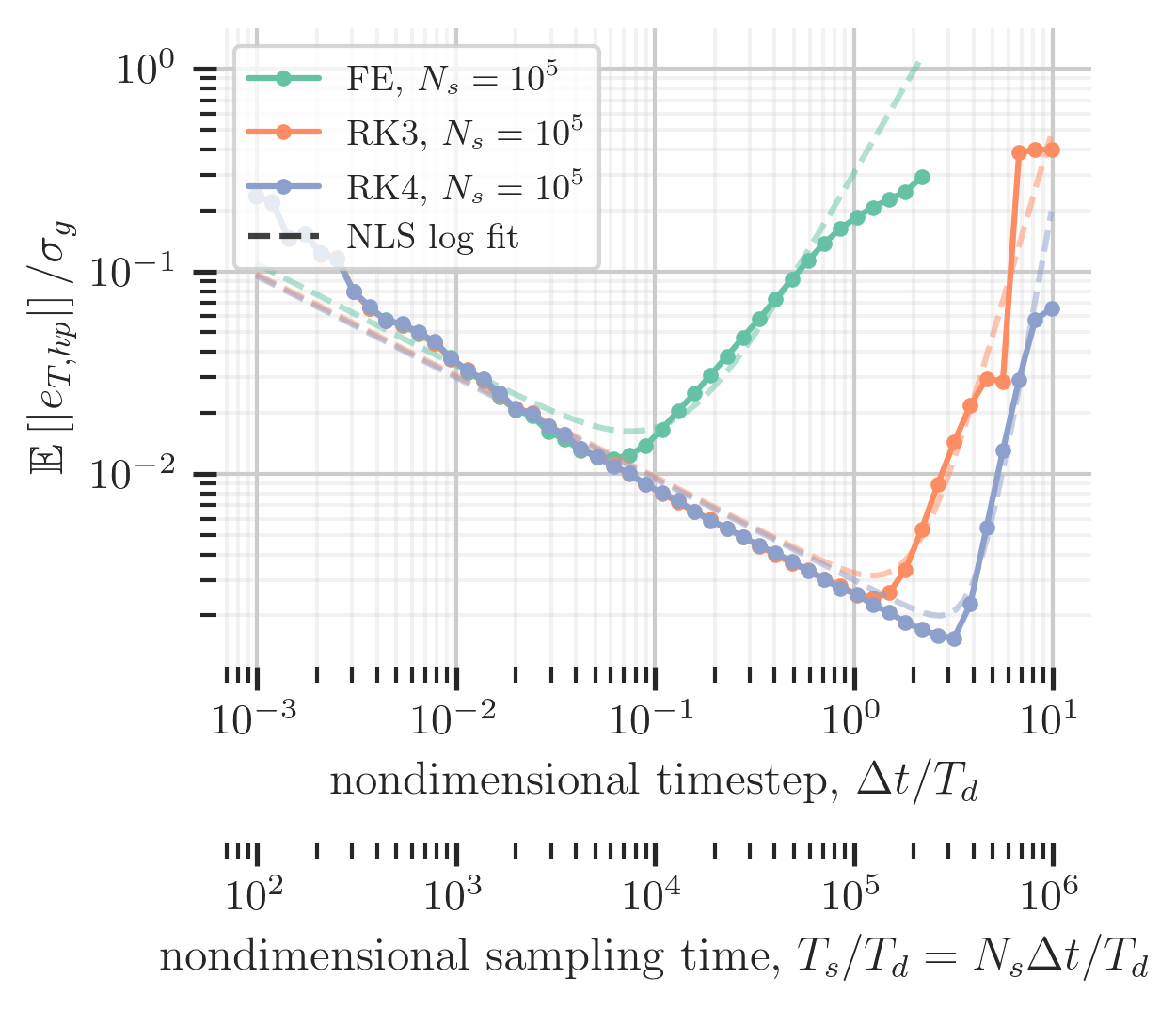}%
    \caption{Expected non-dimensional error as a function of non-dimensional timestep for discretizations of the Lorenz equations. $r= 1/2$ assumed.}%
    \label{fig:errormodel_trade_N100000_nondim}%
\end{figure}

Using the reference simulation, we can also find:
\begin{equation}
  \begin{aligned}
    \variance[J_T] \approx \variance[\JTh] &= 1.1692 \times 10^{-4} \\
    \sigma_g^2 \approx \hat{\sigma}_g^2 &= 74.34804 \pm 0.00018 ,
  \end{aligned}
\end{equation}
where $\hat{\sigma}_g$ is an estimate of the standard deviation of $g$. Together, these allow us to estimate:
\begin{equation}
  \begin{aligned}
    T_d & \approx 1.0170 \times 10^{-2} \\
    \sigma_g & \approx 8.6225 .
  \end{aligned}
  \label{eq:unit_scales}
\end{equation}
With these values, we can plot the non-dimensional error model with fixed $r= 1/2$, which is given for $N_s= 10^5$ in Figure~\ref{fig:errormodel_trade_N100000_nondim}.

We now consider the implications of these results for increasing $N_s$. To focus solely on control of the discretization error, increases in $N_s$ can be used to refine $\dt= T_s/N_s$, with $T_s$ fixed. On the other hand, to focus solely on controlling sampling error, $T_s= N_s \dt$ can be increased, holding $\dt$ fixed.
\begin{figure}[h]
    \centering
    \includegraphics[width= \imgwidth]{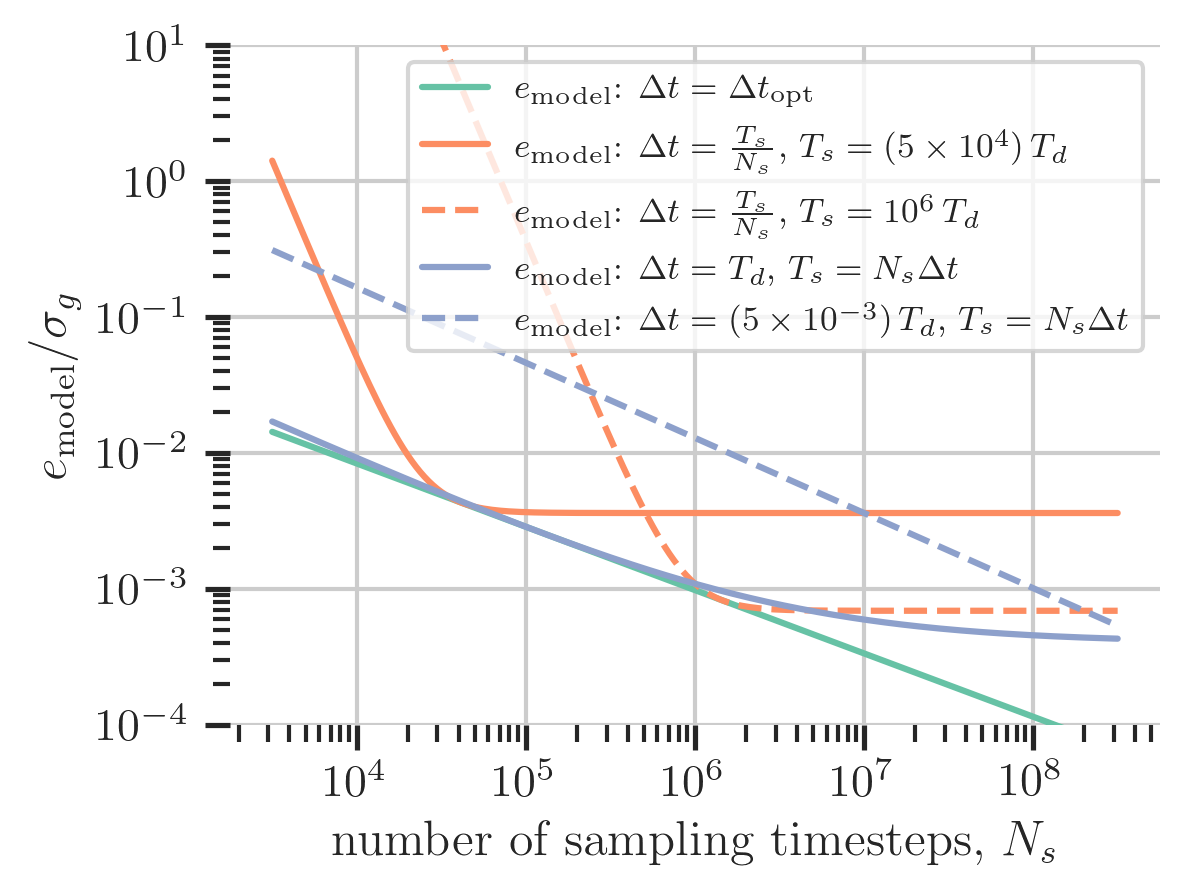}
    \caption{Refinement study comparison for fixed $\dt$, fixed $T_s$, and optimized $\dt$ \& $T_s$ using RK3 discretization to compute the expectation of the Lorenz system output $g= u_2$.}
    \label{fig:refinement_RK3}
\end{figure}
In Figure~\ref{fig:refinement_RK3}, the two approaches are compared with the optimal use of resources. In orange is the discretization error control strategy. In this approach, the simulations converge at a high-order rate in $N_s$ towards the optimal error behavior; once the error reaches this optimum, however, it asymptotes to a constant: statistical errors limit the estimation of $\JTh$. On the other hand, the sampling error control approach is shown in blue. In this approach, the central limit convergence rate of $1/2$ is initially achieved until the error asymptotes to a constant: discretization errors limit the estimation of $\JTh$. In the literature for large simulations, discussed in the introduction, simulations tend to be planned using either the discretization or statistical error control approach. What \eqref{eq:error_model} implies and Figure~\ref{fig:refinement_RK3} demonstrates is that, in fact, there is a particular optimal scheme in which $\dt$ and $T_s$ are simultaneously varied that will extract the most accurate estimate of $\Jinf$ as $N_s$ increases.


\section{Investigation of global discretization error model}
\label{sec:sampling_convergence}

In this section, we show that our simulations of chaotic, ergodic ODEs are consistent with a bounded relationship between the local and global discretization errors. Consider an estimate of the global error based on $N_s$ timesteps:
\begin{equation}
  \hp{e} \approx \frac{1}{N_s} \sum_{n= 0}^{N_s} \sum_{\eta= n}^{N_s} \, \mathcal{G}(t_\eta, t_n) \circ \mathbf{e}_{\mathrm{LT}, p}^{(n)} ,
  \label{eq:error_dynamics}
\end{equation}
where
\begin{equation}
    \mathbf{e}_{\mathrm{LT, p}}^{(n)} \equiv \uhp(t_{n+1}) - \u_{\star}(t_{n+1})
\end{equation}
and $\u_{\star}(t_{n+1})$ is exact solution integrated from $\uhp(t_n)$ through $\dt$:
\begin{equation}
    \u_{\star}(t_{n+1})= \uhp(t_n) + \int_{t_n}^{t_{n+1}} f(\u_{\star}(t)) \diff t .
\end{equation}
In \eqref{eq:error_dynamics}, we have assumed that the error from any given local state perturbation is propagated forward in time by the dynamics, before being transformed into an error in the output; this process is captured by an operator $\mathcal{G}$. Because the effect of local error propagates forward and not backward in time, $\mathcal{G}(t, t_n)= 0$ for $t < t_n$, and moreover we assume that due to ergodicity $\mathcal{G}(t, t_n)= 0$ when $t - t_n \gtrsim T_d$, where $T_d$ is the decorrelation time associated with the attractor. This allows us to write:
\begin{equation}
  \hp{e} \approx \frac{1}{N_s} \sum_{n= 0}^{N_s} \sum_{\eta= n}^{n + T_d/\dt} \, \mathcal{G}(t_\eta, t_n) \circ \mathbf{e}_{\mathrm{LT}, p}^{(n)} .
\end{equation}

Now, we assume that a constant $\mathcal{G}_{\mathrm{max}}$ exists such that:
\begin{equation}
    \abs{\mathcal{G}(t_{\eta}, t_n) \circ \vs} \leq \mathcal{G}_{\mathrm{max}} \norm{\vs}_{\infty} ,
\end{equation}
for all $t_n, t_\eta \in \realno$ and $\vs \in \mathrm{B}(\u(t_n)) \subset \realno[d]$ where $\mathrm{B}(\u)$ is the set of states possible by perturbation of $\u$ that remain in the basin of attraction of the attractor $\attractor$ of $f$. When this is the case, we can create a bound on the magnitude of $\hp{e}$:
\begin{equation}
  \begin{aligned}
    \abs{\hp{e}}
    & \leq \frac{T_d}{\dt} \mathcal{G}_{\mathrm{max}} \frac{1}{N_s} \sum_{n=0}^{N_s} \norm{\mathbf{e}_{\mathrm{LT}, p}^{(n)}}_{\infty} \\
    & \leq \frac{T_d}{\dt} \mathcal{G}_{\mathrm{max}} \max_n \norm{\mathbf{e}_{\mathrm{LT}, p}^{(n)}}_{\infty}
  \end{aligned}
  \label{eq:LTE2global}
\end{equation}

\begin{figure}[h]
    \centering
    \includegraphics[width= \imgwidth]{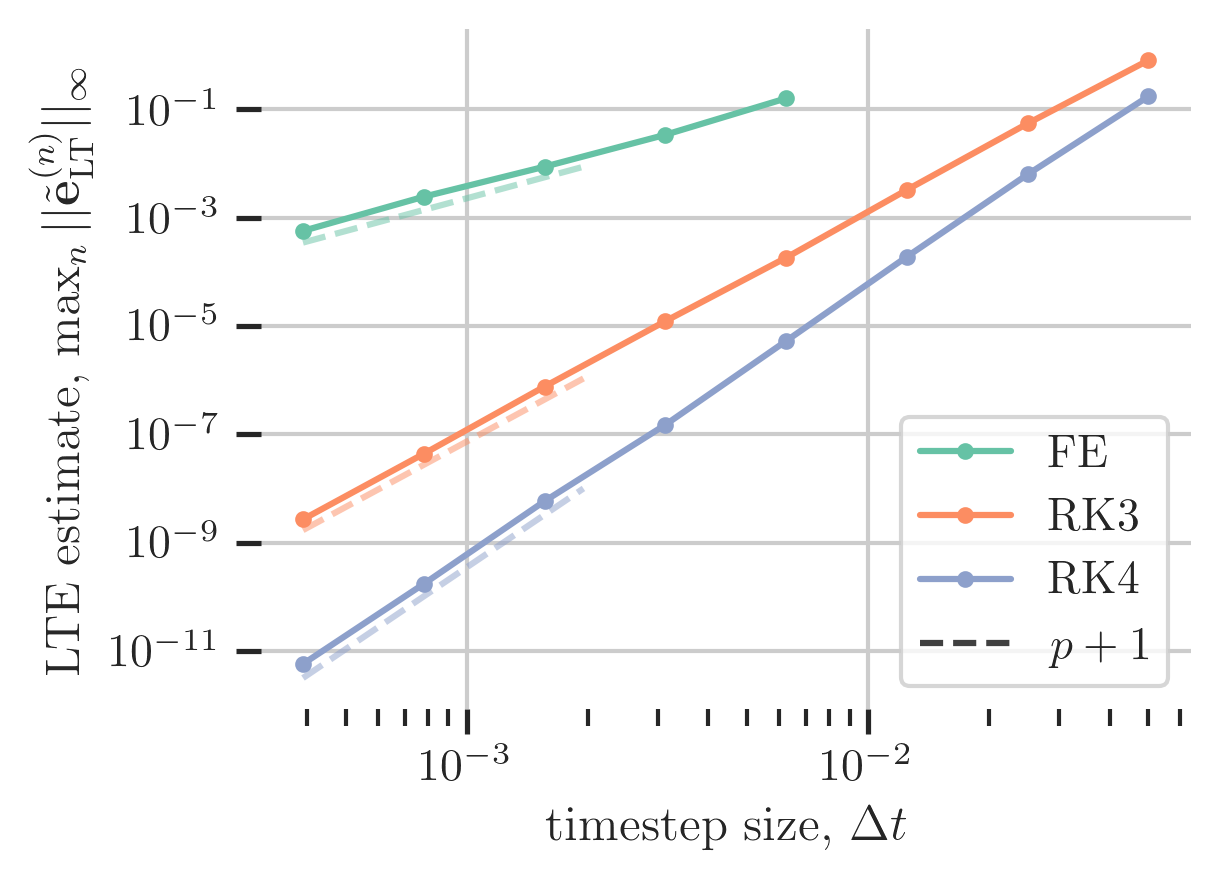}
    \caption{Convergence of estimated local truncation error with respect to $\dt$. Fits to $c_p \dt^{p+1}$ shown (with offset for presentation).}
    \label{fig:LTE_conv}
\end{figure}

We now attempt to bound the value of $\mathcal{G}_{\mathrm{max}}$ for the Lorenz system by approximating the local truncation error. To make an estimate, we compute both the solution at the next timestep as well as a surrogate for the true solution at each timestep: $\uhp(t_{n+1})$ and $\tilde{\u}_{\star}(t_{n+1})$, where the former is computed with one timestep of the method of interest and the latter is always computed with the highest available accuracy method, RK4, and subdividing $t \in [t_{n}, t_{n+1}]$ into ten consecutive timesteps rather than one. Both $\uhp(t_{n+1})$ and $\tilde{\u}_{\star}(t_{n+1})$ are always advanced from $\uhp(t_{n})$. This allows us to estimate $\mathbf{e}_{\mathrm{LT}, p}^{(n+1)}$ locally:
\begin{equation}
    \mathbf{e}_{\mathrm{LT}, p}^{(n+1)} \approx \tilde{\mathbf{e}}_{\mathrm{LT}, p}^{(n+1)}=
        \uhp(t_{n+1}) - \tilde{\u}_{\star}(t_{n+1}) .
    \label{eq:LTE_estimate}
\end{equation}
In Figure~\ref{fig:LTE_conv} we characterize the convergence of local error estimates. Computations are run with $T_s= 100$ and $t_0= 100$ fixed, varying $\dt$. At each timestep, the local truncation error is estimated by computing \eqref{eq:LTE_estimate}. The figure shows the computed $\max_{n} ||\tilde{\mathbf{e}}_{\mathrm{LT}, p}^{(n)}||_{\infty}$ and demonstrates that the expected rate of $(p + 1)$ is nearly exactly achieved.

Using \eqref{eq:LTE2global} we can estimate a bounding value for $\mathcal{G}_{\mathrm{max}}$ by
\begin{equation}
    \mathcal{G}_{\mathrm{max}} \geq \frac{\expect[\abs{\hp{e}}]}{\max_n \norm{\mathbf{e}_{\mathrm{LT}, p}^{(n)}}_{\infty}} \frac{\dt}{T_d}
    = \frac{C_q \dt^{q}}{c_p \dt^{p+1}} \frac{\dt}{T_d} ,
    \label{eq:opG_estimation}
\end{equation}
where $c_p$ is the leading truncation error coefficient fit in Figure~\ref{fig:LTE_conv}, and $C_q$ and $q$ are taken from Table~\ref{tab:NLS_fits_big}. Of course when $q > q_{\mathrm{theory}}= p$, there will be $\dt$ dependence\footnote{In general, we expect $q= p$, but due to cancellation of local errors, $q > p$ occurs in practice for the Lorenz system. In the expected case of $q= p$, we should expect $\mathcal{G}_{\mathrm{max}}= C_q/(c_p T_d)$.}.
However, as \eqref{eq:opG_estimation} requires that the discretization error has an asymptotic behavior, we will only consider $\dt$ in the asymptotic convergence regions given in Table~\ref{tab:fit_bounds} to compute $\mathcal{G}_{\mathrm{max}}$.
\begin{figure}[h]
    \centering%
    \includegraphics[width= \imgwidth]{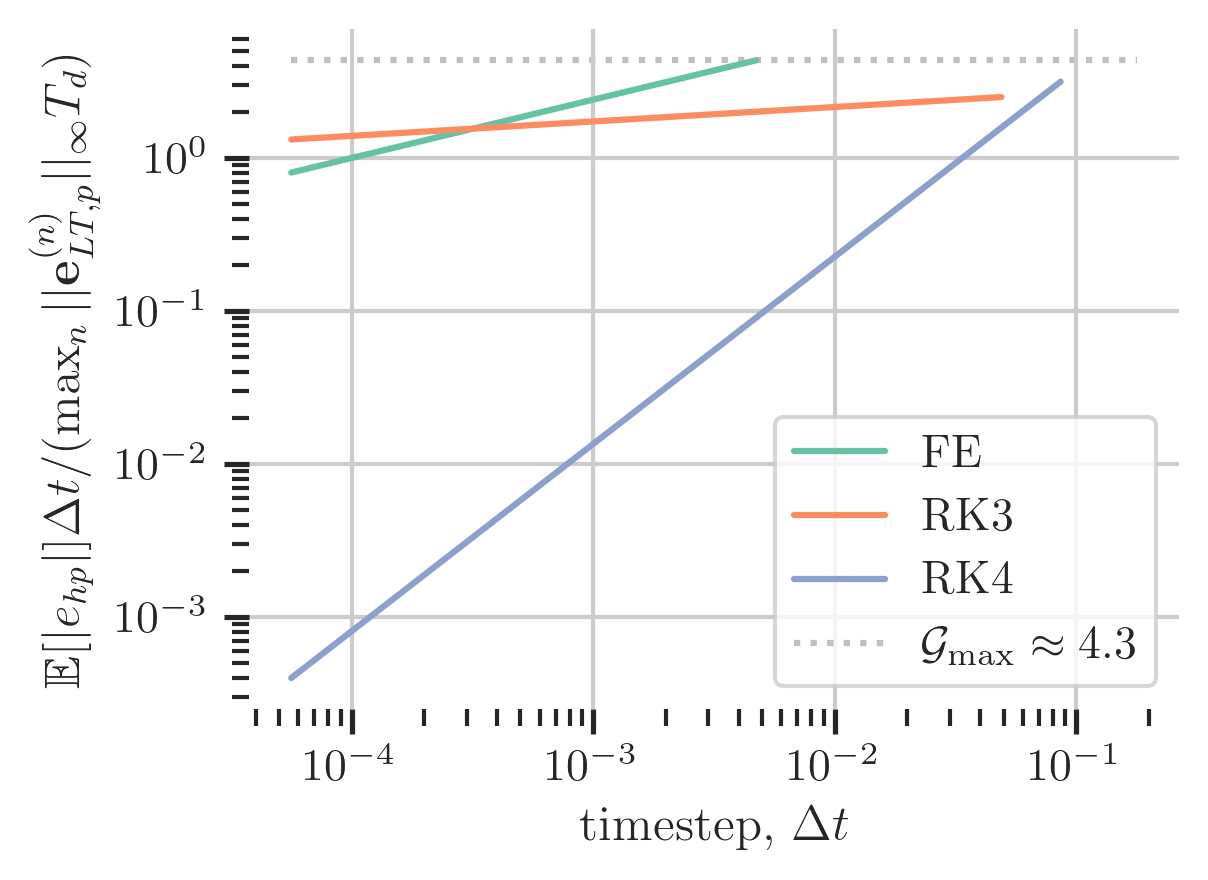}%
    \caption{Estimation of bounding value $\mathcal{G}_{\mathrm{max}}$.}
    \label{fig:opG_estimation}
\end{figure}
In Figure~\ref{fig:opG_estimation}, we show the values of the right-hand side quantity in \eqref{eq:opG_estimation}, which allow us to make an estimate:
\begin{equation}
    \mathcal{G}_{\mathrm{max}} \approx 4.3 .
    \label{eq:opG_est}
\end{equation}

\begin{figure}[h]
    \centering
    \includegraphics[width= \imgwidth]{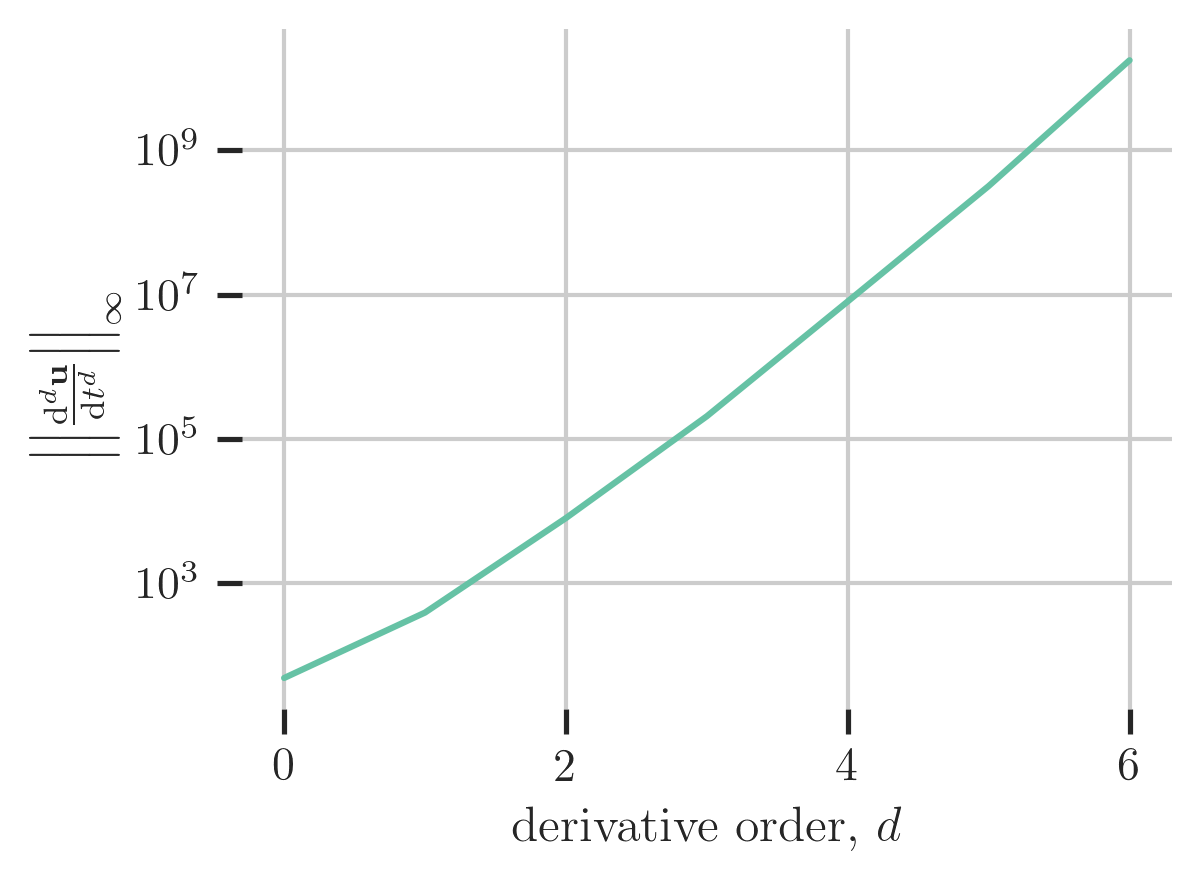}%
    \caption{Norm of analytic derivatives of $\u$ computed on the attractor of $f$. State $\u$ computed with RK4 at $\dt= 10^{-4}$ and $T_s= 1000$ after discarding $t_0= 100$.}
    \label{fig:analytic_derivatives}
\end{figure}

Next, we use classical truncation error estimates \cite{hairer1993solving} to relate the discretization error to properties of the solution. We will assume that the local truncation error is bounded by a form:
\begin{equation}
  \max_n \norm{\mathbf{e}_{\mathrm{LT}, p}^{(n)}}_{\infty} \leq \frac{C_{\mathrm{LT}}}{(p + 1)!} \norm{\dv[p+1]{\u}{t}}_{\infty} \dt^{p+1}
  \label{eq:LTE_theory}
\end{equation}
where $C_{\mathrm{LT}}$ is a local truncation constant term dependent on the numerical method and the $\norm{\cdot}_{\infty}$ in this context refers to the maximum value in time of the inf-norm of a vector-valued, time-dependent quantity $(\cdot)$. The derivatives of $\u(t)$ can be computed by evaluating $f(\u)$ and its derivatives\footnote{Derivatives of $f$ are computed analytically using the chain rule.} using solutions from a reference RK4 solution of the Lorenz system with $T_s= 1000$, $t_0= 100$, and $\dt= 10^{-4}$. Norms of the derivatives are shown in Figure~\ref{fig:analytic_derivatives}.
The resulting values of $C_{\mathrm{LT}}$ that can now be derived by fitting the asymptotic behavior in Figure~\ref{fig:LTE_conv} can be found in Table~\ref{tab:LTE_fit}.
The result of these estimates is that we can reliably bound the global error of a dynamical system as an accumulation of the local errors over a region of correlation.
\begin{table}[h]
    \centering
    \begin{tabular}{c|c|c|c}
        $p$ & $\substack{\text{rate} \\ \text{(observed)}}$ & $C_{\mathrm{LT}} \norm{\dv[p+1]{\u}{t}}_{\infty}$ & $C_{\mathrm{LT}}$ \\
        \hline
        $1$ &   $2.00$ &    $\num{7.61e+03}$ &  $\num{7.33}$ \\
        $3$ &   $4.02$ &    $\num{3.28e+06}$ &  $\num{156}$ \\
        $4$ &   $4.93$ &    $\num{4.50e+07}$ &  $\num{76.5}$
    \end{tabular}
    \caption{Rate and coefficient fit for convergence of local truncation error of discrete Lorenz system. $C_{\mathrm{LT}} \norm{\dv[p+1]{\u}{t}}_{\infty}$ estimated by $c_p (p + 1)!$ using $c_p$ fit from Figure~\ref{fig:LTE_conv}.}
    \label{tab:LTE_fit}
\end{table}

\begin{figure}[h]
  \centering
  \includegraphics[width= \imgwidth]{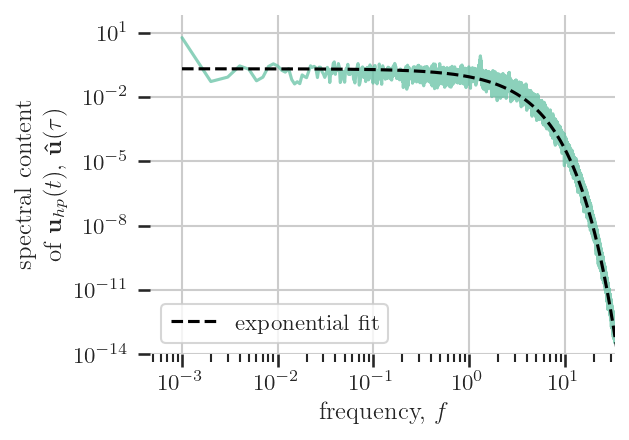}
  \caption{Fourier spectrum of $\u(t)$. Computed with DFT using Hann window function on data from RK4 discretization of Lorenz system with $T_s= 1000$, $t_0= 100$, and $\dt= 10^{-3}$. Gray dashed line: fit assuming $\abs{\hat{\u}(f)} \approx \exp(-a f + b)$ with $a= 0.872$ and $b= 2.58$.}
  \label{fig:spectrum_RK4}
\end{figure}

We now want to consider how the global error behavior demonstrated here might extrapolate to more complicated systems by evaluating the spectral behavior of the Lorenz system. Using a discrete Fourier transform with a Hann window function \cite{harris1978windows}, we perform a spectral analysis on the states of the Lorenz system with a sampling time $T_s= 1000$, $t_0= 100$, and $\dt= 10^{-3}$. The resulting spectrum can be found in Figure~\ref{fig:spectrum_RK4}.
We now want to consider how the demonstrated gloval error behavior demonstrated here results might extrapolate to more complicated systems by evaluating the spectral behavior of the Lorenz system. Using a discrete Fourier transform with a Hann window function \cite{harris1978windows}, we perform a spectral analysis on the states of the Lorenz system with a sampling time $T_s= 1000$, $t_0= 100$, and $\dt= 10^{-3}$. The resulting spectrum can be found in Figure~\ref{fig:spectrum_RK4}.
The Lorenz system tends to have the most content in the frequencies with $f \lesssim 10^{1}$, with a region of exponential decay in the range $1 \lesssim f \lesssim 300$. On scales with $f \gtrsim 300$, machine precision plateaus are observed and omitted here.

The fact that the Lorenz spectrum is an exponentially decreasing function of frequency $f$ makes the use of high-order methods theoretically appealing for the spectral convergence of $hp$-refinement strategies \cite{karniadakis2005spectral}.
Unfortunately, the effect of statistical error in \eqref{eq:error_optimal_nondim} limits the impact of this exponential decay, such that the benefits of higher-order discretization methods are limited compared to their steady-state and non-chaotic application.
\begin{figure}[h]
    \centering
    \includegraphics[width= \imgwidth]{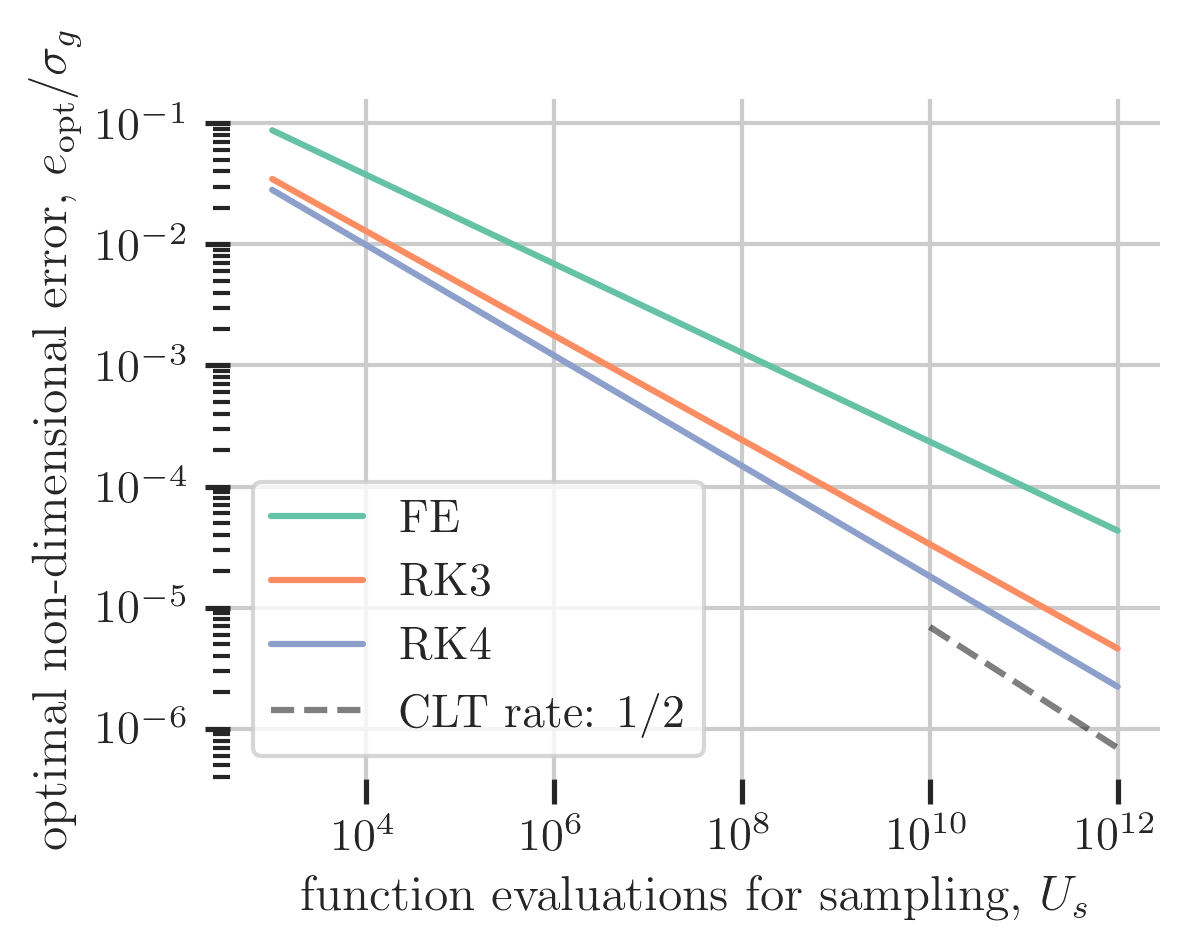}
    \caption{Convergence of optimal error with sampling costs for FE, RK3, and RK4 discretizations of the Lorenz output $g= u_2$. Asymptotic $-1/2$ rate implied by central limit theorem shown.}
    \label{fig:sampling_convergence}
\end{figure}
The convergence to the central limit rates can be seen in Figure~\ref{fig:sampling_convergence}, which shows the convergence of \eqref{eq:error_optimal_nondim} with the total sampling cost. The effect of increasing order improves the convergence rate in \eqref{eq:error_optimal_nondim} towards the CLT-implied asymptotic rate of $-1/2$, as well as decreasing the value of the leading constant and the error never achieves the spectral rates possible with $hp$-refinement in the steady case. Nevertheless, the cost to achieve a given amount of error in expectation-- in terms of function evaluations-- is significantly less with higher-order methods.
Managing to achieve 1\% non-dimensional error in expectation is possible with RK4 at a cost ten times less than would be possible using FE; that factor grows larger than 100 when the tolerance is tightened to $10^{-4}$.


\section{Impact of ensemble averaging and spin-up}

In this section, we will consider how the error behaves when ensemble averaging (over multiple parallel instances) and when spin-up effects are present.

\subsection{Ensemble averaging on the attractor}

Sampling error can be reduced at a fixed wall clock time by ensemble averaging across multiple parallel processes \cite{makarashvili2017performance}. Consider a Monte Carlo approach to approximate $\Jinf$ with a set of $\Mens$ independent realizations:
\begin{equation}
    J_\mathrm{MC}= \frac{1}{\Mens} \sum_{m= 1}^{\Mens} \JTh^{(m)} .
    \label{eq:output_mc}
\end{equation}
We can write a modified version of \eqref{eq:error_model} to approximate the error that we expect in the Monte Carlo estimator in \eqref{eq:output_mc}:
\begin{widetext}
\begin{equation}
    \expect[|J_\mathrm{MC} - \Jinf|] \approx
    e_{\mathrm{model, MC}}=
    C_q \dtMC^q + \frac{A_0}{\sqrt{\Mens}} \TsMC^{-r} ,
    \label{eq:error_model_mc}
\end{equation}
\noindent
with an equivalent non-dimensional version, assuming $r \to 1/2$:
\begin{equation}
  \left( \frac{\emodel}{\sigma_g} \right)_{\mathrm{MC}}=
      \frac{C_q T_d^q}{\sigma_g} \left( \frac{\dt}{T_d} \right)^q
      + \sqrt{\frac{2}{\pi}} \Mens^{-\frac{1}{2}} \left( \frac{T_s}{T_d} \right)^{-\frac{1}{2}} ,
  \label{eq:error_model_mc_nondim}
\end{equation}
and an optimum given by
\begin{equation}
  \left( \frac{\emodel}{\sigma_g} \right)_{\mathrm{MC}, \mathrm{opt}}= \\
  \left( \frac{1}{2 \pi} \right)^{\frac{q}{2q+1}} \left( 2 + \frac{1}{q} \right)
      \left( \frac{q C_q T_d^q}{\sigma_g} \right)^{\frac{1}{2q+1}} \Mens^{-\frac{q}{2q+1}}
      N_s^{-\frac{q}{2q+1}} ,
  \label{eq:error_optimal_mc_nondim}
\end{equation}
at
\begin{equation}
  \left( \frac{\dt}{T_d} \right)_{\mathrm{opt}}=
  \left( \frac{1}{2 \pi} \right)^{\frac{1}{2q+1}} \left(
  \frac{q C_q T_d^q}{\sigma_g} \right)^{-\frac{2}{2q+1}} \Mens^{-\frac{1}{2q+1}} N_s^{-\frac{1}{2q+1}} ,
  \label{eq:dt_opt_MC_nondim}
\end{equation}
and
\begin{equation}
  \left( \frac{T_s}{T_d} \right)_{\mathrm{opt}}= \left( \frac{1}{2 \pi}
  \right)^{\frac{1}{2q+1}} \left( \frac{q C_q T_d^q}{\sigma_g}
  \right)^{-\frac{2}{2q+1}} \Mens^{-\frac{1}{2q+1}} N_s^{\frac{2q}{2q+1}} .
  \label{eq:Ts_opt_MC_nondim}
\end{equation}
\end{widetext}

Equation~\ref{eq:error_optimal_mc_nondim} shows that, for finite values of $q$, the Monte Carlo method will have a mitigated return compared to its purely stochastic application as in \citet{makarashvili2017performance}; the optimal error scales as $\Mens^{-q/(2 q + 1)}$ as opposed to $\Mens^{-1/2}$. However, parallelization can achieve perfect scaling in the expected error, in the sense that the effect of running $\Mens$ ensembles with $N_s$ sampling timesteps each will have an equivalent error in expectation to simulating $\Mens N_s$ timesteps in serial. As $\Mens$ is varied on the set of optimal solutions, \eqref{eq:dt_opt_MC_nondim} and \eqref{eq:Ts_opt_MC_nondim} indicate that the timestep and sampling time should be adjusted with the same factor $\Mens^{-1/(2q+1)}$ to achieve perfect scaling.

\subsection{Spin-up transient modeling}

So far, we have considered the error and cost on the attractor, neglecting the impact of ``spin-up'' from $t= 0$ to $t= t_0$. This spin-up is necessary because simulations of ergodic systems invariably need some time for the state to proceed onto the attractor from the initial condition.

Consider $\u(t)$, a solution of the ergodic chaotic system $f$ from an arbitrary initial condition $\u(0)= \uIC$ in the basin of attraction of an attractor, $\attractor$. The existence of the attractor implies the non-linear stability of the system, such that all $\uIC$ will converge to trajectories on the attractor $\attractor$. Denote by $\u^{\attractor}(t)$ a trajectory that is on the attractor for all $t$ and to which $\u(t)$ collapses as $t \to \infty$. The perturbation $\delta \u^{\attractor}(t) \equiv \u(t) - \u^{\attractor}(t)$ that describes the IC, therefore, exists in a stable subspace of perturbations to $\u^{\attractor}$ and can be associated with the negative Lyapunov exponents of the system. Thus, we can assume that such perturbations are governed asymptotically by
\begin{equation}
    \norm{\delta \u^{\attractor}(t)} \lesssim \exp(-\frac{t}{T_{\lambda}}) ,
    \label{eq:delta_lyapunov_transient}
\end{equation}
with $T_{\lambda}$ a characteristic time associated with the stable Lyapunov modes. In practice, we are interested in averages of quantities on the attractor $g(\u^{\attractor}(t))$, but we can only calculate quantities $g(\u(t))$, that will include some effect-- if small-- of the spin-up transient.

Next, we seek to quantify the effect of this gap on estimates $J_T \approx \Jinf$. Consider the computation of $J_T$. In \eqref{eq:error_total}, we have effectively found an estimate of
\begin{equation}
    J_T^{\attractor}= \int_{t_0}^{t_0 + T_s} g(\u^{\attractor}(t)) \diff t ,
\end{equation}
by choosing $t_0$ sufficiently large. We now want to consider an error model of the form:
\begin{equation}
    e_{T, hp}= \underbrace{(\JTh - J_T)}_{\hp{e}}
    + \underbrace{(J_T - J_T^{\attractor})}_{e_{\lambda}}
    + \underbrace{(J_T^{\attractor} - \Jinf)}_{e_T}
    \label{eq:error_total_spin-up}
\end{equation}
where a new error $e_{\lambda}$ is introduced, associated with the spin-up transient. The model for $e_T$ in \eqref{eq:error_stat} will apply without modification, while the model for $\hp{e}$ will be subject to slightly different assumptions. Where in \eqref{eq:error_disc}, $C_q$ was bounded by the value on the attractor, $\attractor$, here we must assume that $C_q$ is bounded from $t= 0$ to $t= t_0 + T_s$, including both the attractor \emph{and} the transient part of the trajectory. We only require that the transient part be in the basin of attraction of $\attractor$, $B(\attractor)$. We assume that a model of the form used in \eqref{eq:error_disc} applies in expectation when the transient component is included.

Next, we concentrate on $e_{\lambda}$:
\begin{equation}
    J_T - J_T^{\attractor}= \int_{t_0}^{t_0 + T_s} \left( g(\u(t)) - g(\u^{\attractor}(t)) \right) \diff t .
\end{equation}
We now assume that, like $\u$, $g$ will decay exponentially in $t$ as \eqref{eq:delta_lyapunov_transient}, such that
\begin{equation}
    g(\u(t)) - g(\u^{\attractor}(t)) \equiv \delta g^{\attractor}(t) \approx A_{\lambda} \exp(-\frac{t}{T_{\lambda}})
    \label{eq:output_state_diff}
\end{equation}
will apply for $t \in [0, \infty)$, with $A_{\lambda}$ a constant that can be related to the deviation between $g(\u(0))$ and $g(\u^{\attractor}(0))$.

From this assumption,
\begin{equation}
    \begin{aligned}
        e_{\lambda} &= \frac{1}{T_s} \int_{t_0}^{t_0 + T_s} g(\u(t)) - g(\u^{\attractor}(t)) \diff t \\
        & \approx \frac{1}{T_s} \int_{t_0}^{t_0 + T_s} A_{\lambda} \exp(-\frac{t}{T_{\lambda}}) \diff t \\
        &= A_{\lambda} \frac{T_{\lambda}}{T_s} \exp(-\frac{t_0}{T_{\lambda}}) \left( 1 - \exp(-\frac{T_s}{T_{\lambda}}) \right)
    \end{aligned}
\end{equation}
Taking the absolute value, we can find a bounding model:
\begin{equation}
    |e_{\lambda}|= |A_{\lambda}| \frac{T_{\lambda}}{T_s} \exp(-\frac{t_0}{T_{\lambda}}) .
    \label{eq:abs_transient_model}
\end{equation}
As before, manipulation of \eqref{eq:error_total_spin-up} allows
\begin{align}
    |e_{T, hp}| & = |\hp{e} + e_{\lambda} + e_T| \\
    & \leq |\hp{e}| + |e_{\lambda}| + |e_T| .
    \label{eq:error_total_spin-up_triangle}
\end{align}
Now, we take an expectation of the absolute value of $e_{T, hp}$:
\begin{equation}
    \expect[|e_{T, hp}|] \leq \expect_{B(\attractor)}[|\hp{e}|] + \expect_{\mathrm{IC}}[|e_{\lambda}|] + \expect_{\attractor}[|e_T|] ,
    \label{eq:error_total_spin-up_triangle_expect}
\end{equation}
where $\expect_{B(\attractor)}$ gives the expectation on the basin of attraction of $\attractor$. Here, the expectation of $\abs{\eTh}$ doesn't reduce to an expectation \emph{on the attractor}. The statistical term is handled on the attractor as before, and we have assumed that the discretization error is bounded by the same form in expectation on $B(\attractor)$ as on $\attractor$. Finally, the expectation of $\abs{e_{\lambda}}$ is taken on the set of initial conditions used. This allows us to take the expectation of \eqref{eq:abs_transient_model} to complete \eqref{eq:error_total_spin-up_triangle_expect}. Because we anticipate $T_{\lambda}$ will be bounded by a constant for a given system, this is given by:
\begin{equation}
    \expect_{\mathrm{IC}}[|e_{\lambda}|]= \expect_{\mathrm{IC}}[|A_{\lambda}|] \frac{T_{\lambda}}{T_s} \exp(-\frac{t_0}{T_{\lambda}}) .
    \label{eq:exp_transient_model}
\end{equation}
If a $A_{\lambda}$ and $T_{\lambda}$ can be identified by observation of $g(\u(t))$ given an initial condition $\u_{\mathrm{IC}}$, $\abs{e_{\lambda}}$ is no longer stochastic and the $\expect[|e_{T, hp}|] \to |e_{T, hp}|$ as in \eqref{eq:abs_transient_model}.

\begin{widetext}
Putting all the pieces together, we can now give an error model that incorporates the effects of spin-up and ensemble estimation:
\begin{equation}
    e_{\mathrm{model, MC}}= \tilde{A}_{\lambda} \frac{T_{\lambda}}{\TsMC} \exp(-\frac{t_0}{T_{\lambda}}) + C_q \dtMC^q + \frac{A_0}{\sqrt{\Mens}} \TsMC^{-r} ,
    \label{eq:error_model_mc_spin-up}
\end{equation}
where $\tilde{A}_{\lambda}$ can be either estimated on an instance-by-instance basis or by estimating the expectation on the family of initial conditions.
\end{widetext}
Under this model, $e_{\lambda}$ will scale with the exponent of a large negative value when $t_0 \gg T_{\lambda}$. Even when $t_0 \not\gg T_{\lambda}$, \eqref{eq:abs_transient_model} suggests that the decay-induced error term will still scale with $T_s^{-1}$, faster than the expected CLT rate of $T_s^{-1/2}$, and thus it will be dominated it as $T_s \gg 1$. This implies two ``paths'' to controlling spin-up errors: either choosing $t_0$ long enough to shrink the mean offset error from $t_0$, or choosing $T_s$ long enough so that the mean offset contribution to the simulation error is small in spite of the error at $t= t_0$.

\subsection{Identification of spin-up transient model}

We will now develop a method to fit the error model. In order to do so, consider observations $g_n \equiv g(t_n)$ and $g_n^{\attractor} \equiv g^{\attractor}(t_n)$ for $t_n$ in $\{t_0, t_0 + N_{\mathrm{skip}} \dt, \ldots, t_0 + T_s\}$. We will assume that $N_{\mathrm{skip}}$ is large enough that the solution at each $t_n$ is effectively independent. If this is the case, then we can assume that each $g_n^{\attractor}$ will be an independent and identically distributed (i.i.d.) draw from a bounded, stationary distribution with mean $\Jinf$. The distributions of $g^{\attractor}(\u(t))$ and $g(\u(t))$, in general, are not known. In order to facilitate an estimate of the mean behavior, we will assume $g_n^{\attractor}$ are i.i.d. draws from a normal distribution with mean value $\Jinf$. Then, we have:
\begin{equation}
    g_n \sim \normaldist(\Jinf + \delta g^{\attractor}(t_n), \sigma_g^2) ,
    \label{eq:mean_offset_rv}
\end{equation}
where the relationship between $g_n$ and $g_n^{\attractor}$ is taken from \eqref{eq:output_state_diff}.

In order to understand the implications of this model, we can use set of reference RK4 simulations of the Lorenz system with $N_t= 10^5$ timesteps sampled without spin-up over a period $T= 100$ from initial conditions similar to those given in \eqref{eq:random_ic}, with a scaled-up standard deviation of 100 in all three variables to highlight the initial transient. In order to treat each of $\Jinf$, $\sigma_g$, $A_{\lambda}$, and $T_{\lambda}$ in \eqref{eq:mean_offset_rv} as unknowns, we use Hamiltonian Monte Carlo with the likelihood function implied by \eqref{eq:mean_offset_rv}. We discard from $t= 0$ to $t= 5$, then take \num{10000} equispaced samples from $t= 5$ to $t= 100$. For prior models, we start by computing na\"{i}ve estimators of the mean and standard deviation of the trace, $\tilde{J}$ and $\tilde{\sigma}$ using the downsampled trace signal $\{ g_n \}$, then use:
\begin{equation}
    \begin{aligned}
        \Jinf &\sim \mathcal{N}(\tilde{J}, \tilde{\sigma}^2) \\
        \sigma_g &\sim \Gamma(\alpha_{\sigma}, \beta_{\sigma}) \\
        A_{\lambda} &\sim \mathcal{N}(0, \max(\hp{g}) - \min(\hp{g})) \\
        T_{\lambda} &\sim \Gamma(\alpha_T, \beta_T)
    \end{aligned}
\end{equation}
where
\begin{equation*}
    \begin{aligned}
        (\alpha_{\sigma}, \beta_{\sigma}) & \impliedby \left( \mu_{\sigma}= \tilde{\sigma}, \sigma_{\sigma}= \frac{\tilde{\sigma}}{10} \right) \\
        (\alpha_T, \beta_T) & \impliedby \left( \mu_T= 10.0, \sigma_T= 10.0 \right) .
    \end{aligned}
\end{equation*}
It should be noted that in this specification, the Bayesian fit only requires a user-supplied prior for the decay time and for the uncertainty in the standard deviation, assumptions upon which the fitting method only requires be reasonable.

\begin{figure}[h]
    \centering
    \includegraphics[width= \imgwidth]{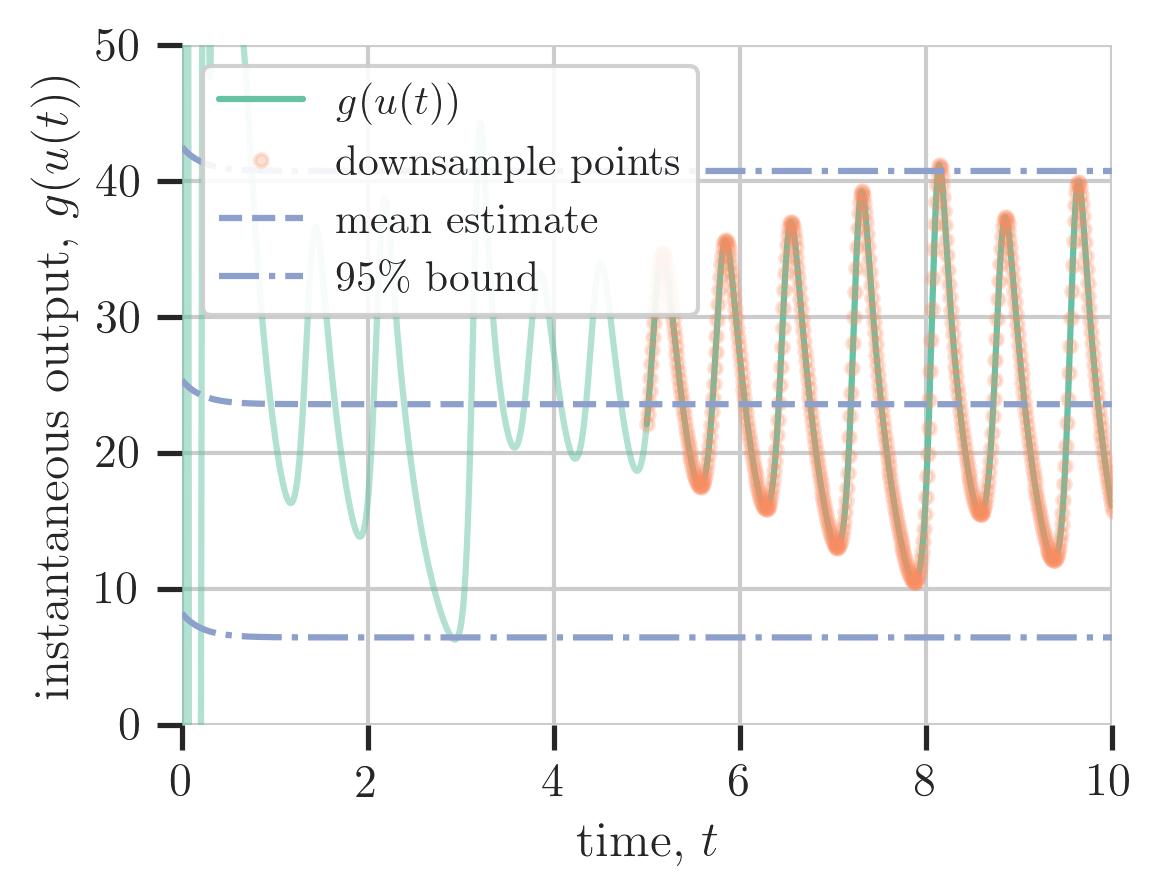}
    \caption{$g= u_2(t)$ trace in transient region, with Bayesian method fit}
    \label{fig:trace_decay_fit_sample}
\end{figure}

\begin{figure}[h]
  \centering
  \includegraphics[width= \imgwidth]{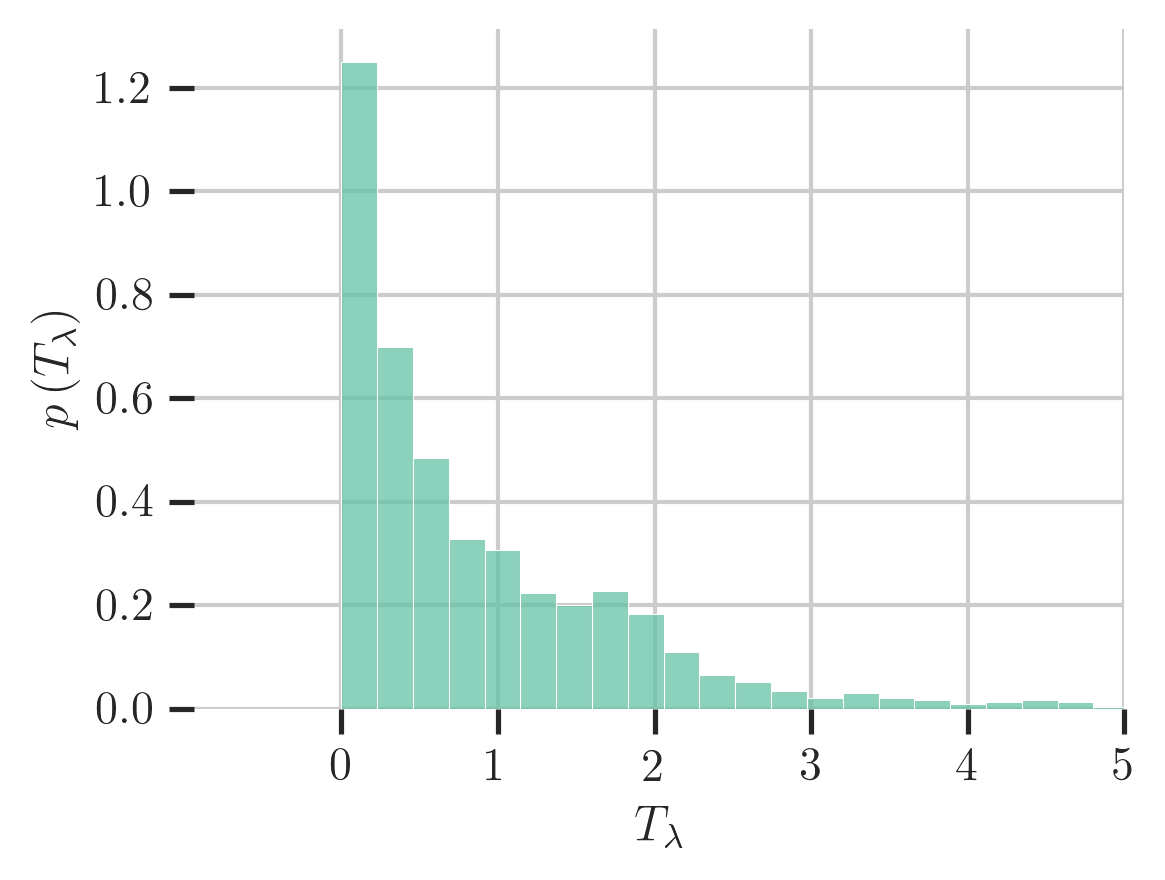}
  \caption{MAP estimate $T_{\lambda}$ for Lorenz system transient. Collected over 1000 Lorenz trajectories with $\dt= 10^{-2}$, $T_s= 100$, and randomized $\uIC$. Outliers truncated, greater than 97\% of data in pictured range.}
  \label{fig:Tlambda_histogram}
\end{figure}

\begin{figure}[h]
  \centering
  \includegraphics[width= \imgwidth]{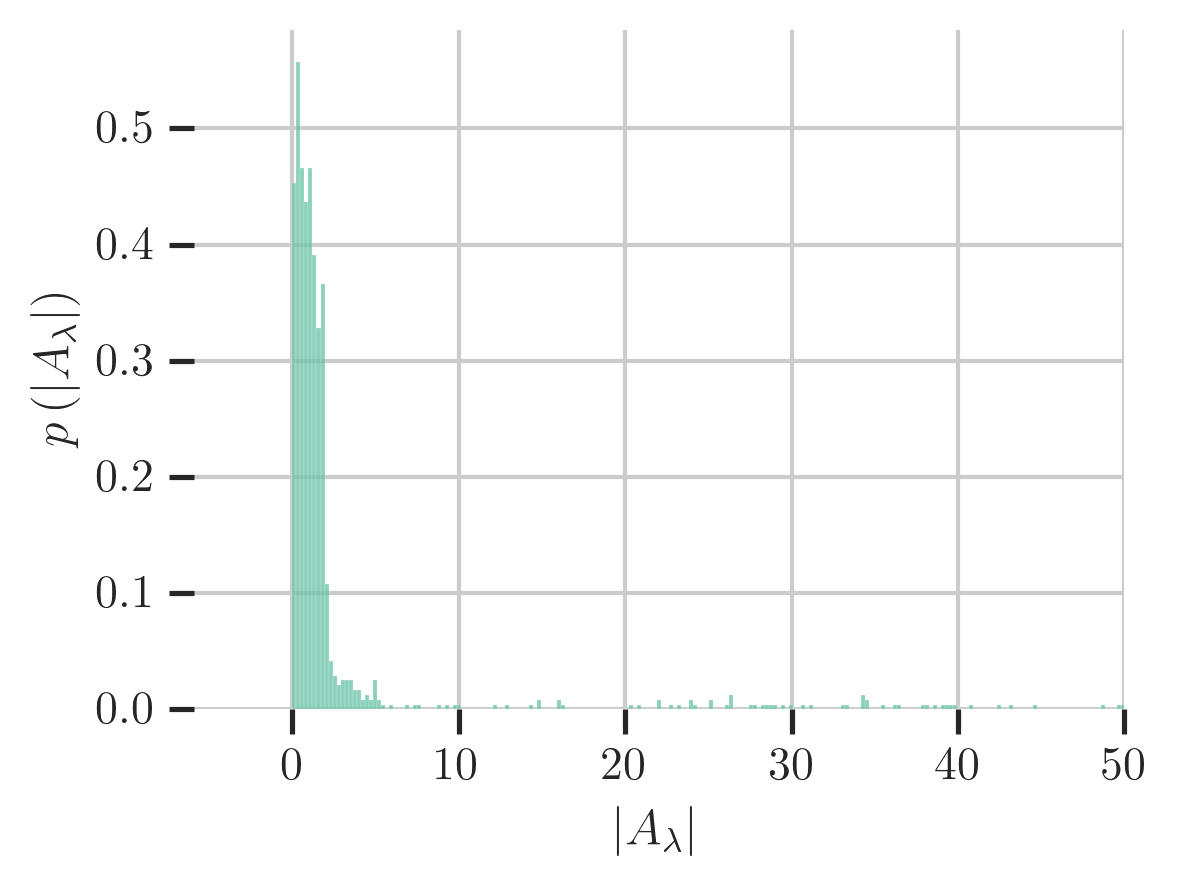}
  \caption{MAP estimate $\abs{A_{\lambda}}$ for Lorenz system transient. Collected over 1000 Lorenz trajectories with $\dt= 10^{-2}$, $T_s= 100$, and randomized $\uIC$. Outliers truncated, greater than 97\% of data in pictured range.}
  \label{fig:absAlambda_histogram}
\end{figure}

A sample fit and trace are found in Figure~\ref{fig:trace_decay_fit_sample}, for which the maximum \apost{} estimate gives $T_{\lambda}= 0.312$ and $A_{\lambda}= -0.925$.
For the Lorenz system, the initial transient onto the attractor is very rapid, almost negligible. Applying the Bayesian fit procedure to an ensemble of 1000 runs generated in the same way as Figure~\ref{fig:trace_decay_fit_sample} we can find maximum \apost{} (MAP) estimates of the variables $T_{\lambda}$ and $\abs{A_{\lambda}}$ in the decay model. In Figures~\ref{fig:Tlambda_histogram} and~\ref{fig:absAlambda_histogram}, histograms of these variables are shown, which are needed to determine \eqref{eq:abs_transient_model}.
We can see that the fit procedure identifies values:
\begin{equation}
  \begin{aligned}
    T_{\lambda} & < 4.03 \\
    |A_{\lambda}| & < 38.7
  \end{aligned}
  \label{eq:transient_model_estimates}
\end{equation}
for greater than 97\% of initial conditions, up to two standard deviations above the mean. Using these values as a conservative estimate for the mean offset, we can now model the effect of the transient behavior.


\section{Optimal time-stepping including spin-up}

Now we can consider how the cost and error impact of spin-up is incorporated into the model for error at a fixed cost. The spin-up time requires the use of $N_0$ timesteps:
\begin{equation}
    N_0= \left\lceil \frac{t_0}{\dtMC} \right\rceil \approx \frac{t_0}{\dtMC} .
\end{equation}
With $N$ the total number of timesteps used, given by:
\begin{equation}
    N= N_0 + \NMC= \frac{t_0}{\dtMC} + \NMC ,
    \label{eq:wallclock_basic}
\end{equation}
where $\NMC$ is the number of timesteps during sampling for $t_0$ to $t_0 + T_s$ on a given instance.

\begin{widetext}
By normalizing \eqref{eq:error_model_mc_spin-up} then substituting \eqref{eq:wallclock_basic}, we arrive at a transient-inclusive non-dimensional model for the error:
\begin{equation}
    \begin{aligned}
        \left( \frac{\emodel}{\sigma_g} \right)_{\mathrm{MC}} &=
                \frac{\tilde{A}_{\lambda}}{\sigma_g} \frac{T_{\lambda}}{T_d} \left( N \left( \frac{\dt}{T_d} \right)_{\mathrm{MC}} - \frac{t_0}{T_d} \right)^{-1} \exp(-\frac{t_0/T_d}{T_{\lambda}/T_d}) \\
        & \quad + \frac{C_q T_d^q}{\sigma_g} \left( \frac{\dt}{T_d} \right)_{\mathrm{MC}}^q
                + \sqrt{\frac{2}{\pi}} \Mens^{-\frac{1}{2}} \left( N \left( \frac{\dt}{T_d} \right)_{\mathrm{MC}} - \frac{t_0}{T_d} \right)^{-\frac{1}{2}} .
    \end{aligned}
    \label{eq:error_model_constrained}
\end{equation}
\end{widetext}
Using this result, we can solve numerically for $\dtMCopt$ and $\eMCopt$ via \eqref{eq:error_model_constrained}.

Consider a Lorenz simulation on which a budget of $U= p N= \num{1.2e6}$ right-hand side evaluations are available on each of $\Mens$ parallel processors. We start by studying the error under \eqref{eq:error_model_constrained} as $\dt$ and $t_0$ vary with a conservative estimate for the transient behavior using the bounding values in \eqref{eq:transient_model_estimates}.
\begin{figure}[h]
    \centering
    \includegraphics[width= \imgwidth]{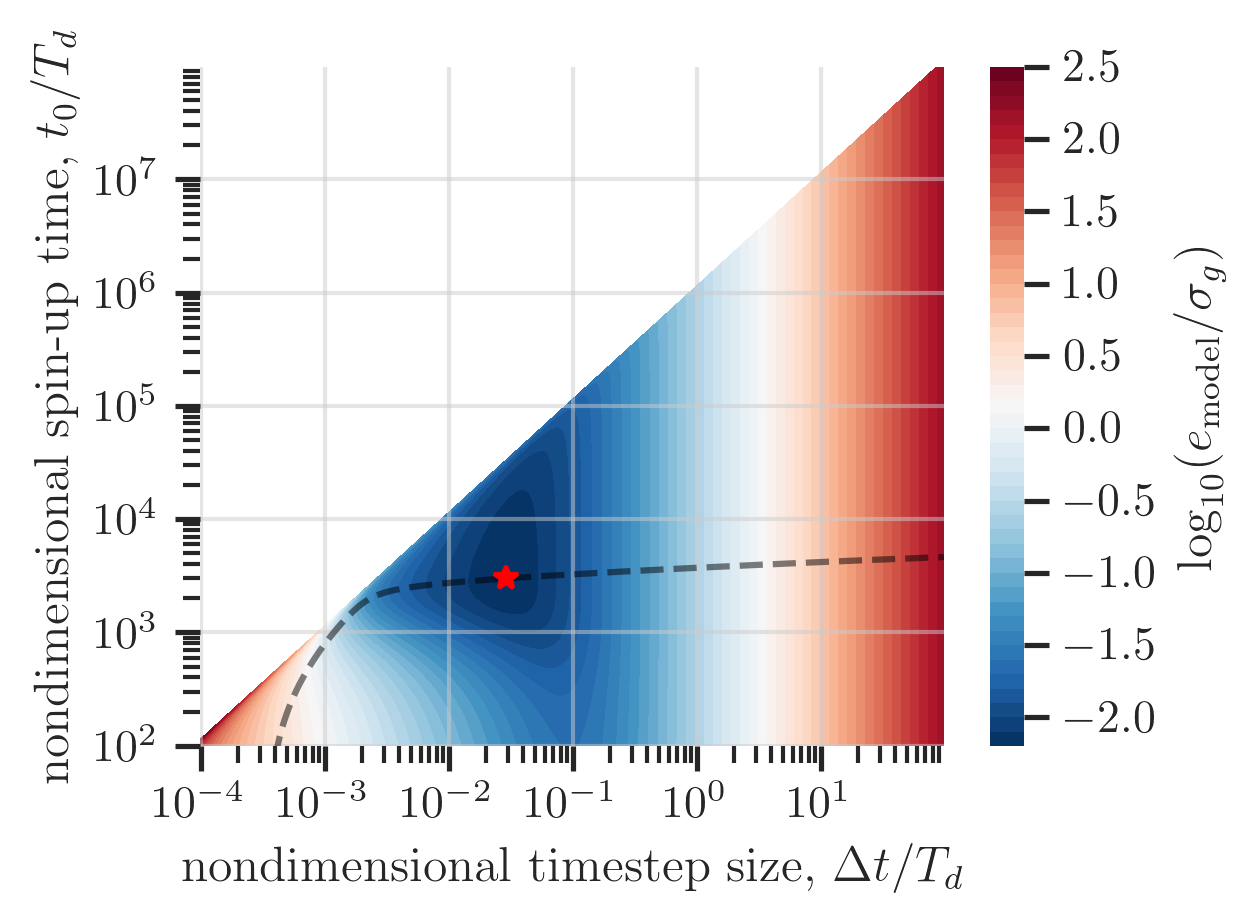}
    \caption{Dependence of normalized error expectation $\eMC/\sigma_g$ on normalized timestep $\dt/T_d$ and normalized spin-up time $t_0/T_d$ with total cost set at $U= \num{1.2e6}$ for Forward Euler. Red star denotes optimum, dashed line indicates optimal $t_0$ given $\dt$.}
    \label{fig:totalcost_dt_t0_emod_FE}
\end{figure}
\begin{figure}[h]
    \centering
    \includegraphics[width= \imgwidth]{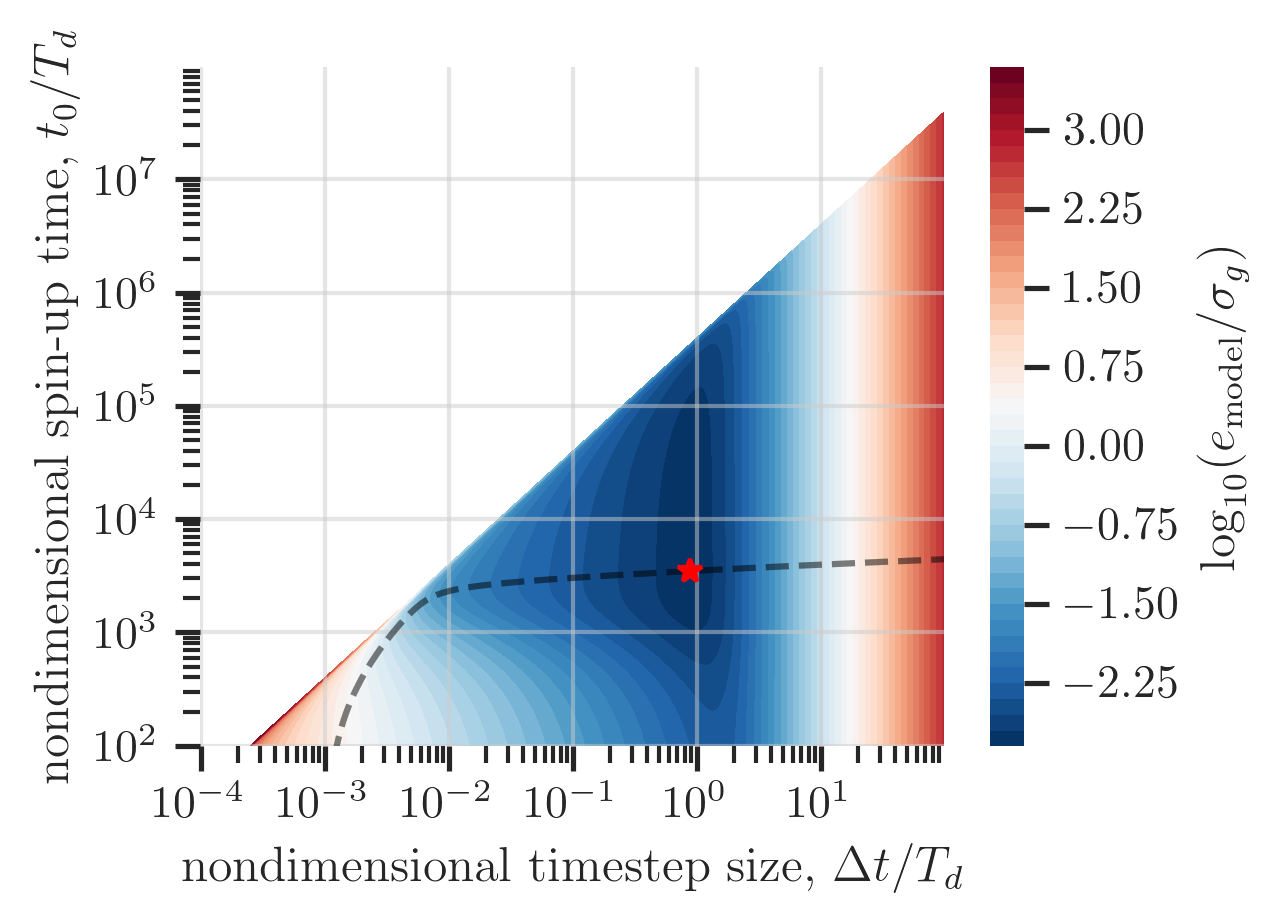}
    \caption{Dependence of normalized error expectation $\eMC/\sigma_g$ on normalized timestep $\dt/T_d$ and normalized spin-up time $t_0/T_d$ with total cost set at $U= \num{1.2e6}$ for \nth{3}-order Runge Kutta. Red star denotes optimum, dashed line indicates optimal $t_0$ given $\dt$.}
    \label{fig:totalcost_dt_t0_emod_RK3}
\end{figure}
\begin{figure}[h]
    \centering
    \includegraphics[width= \imgwidth]{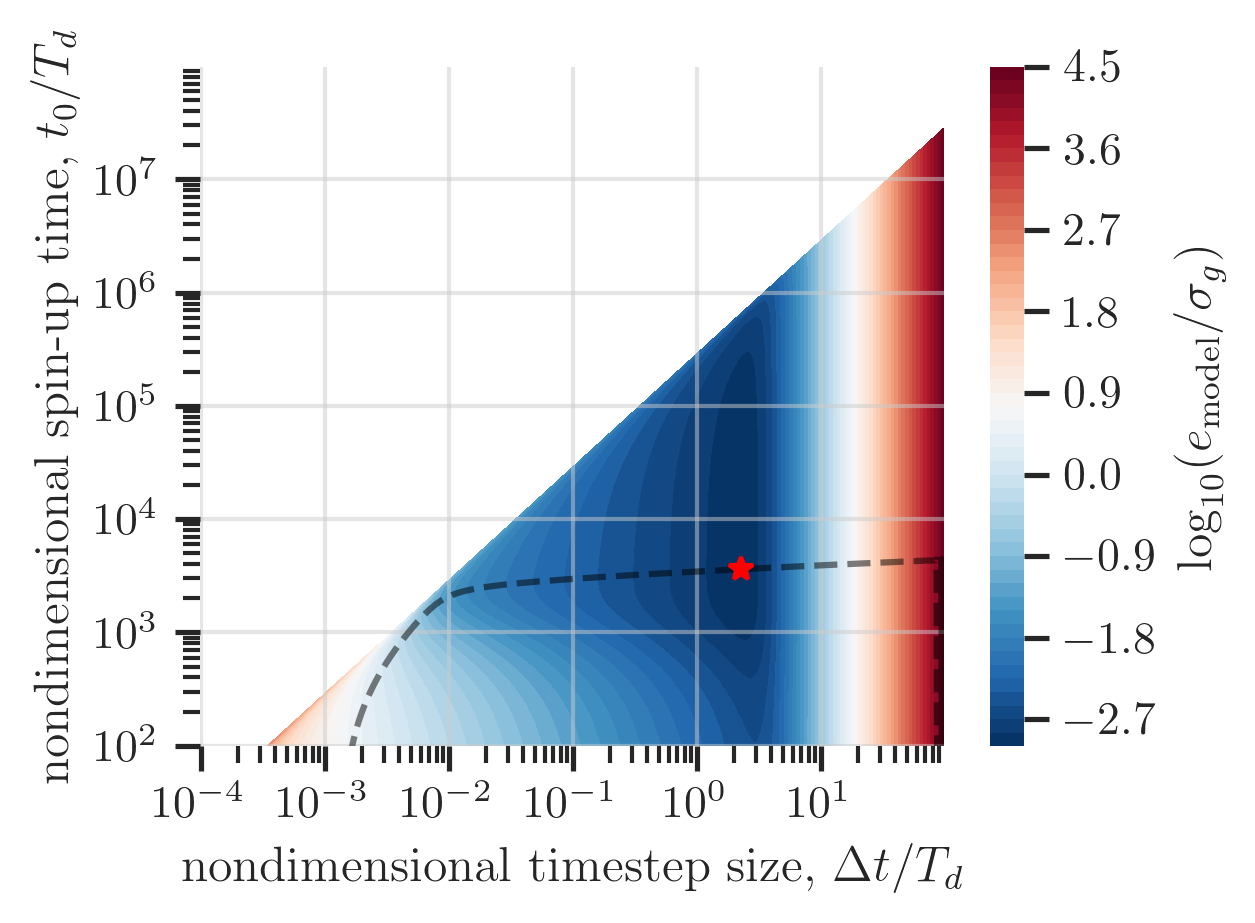}
    \caption{Dependence of normalized error expectation $\eMC/\sigma_g$ on normalized timestep $\dt/T_d$ and normalized spin-up time $t_0/T_d$ with total cost set at $U= \num{1.2e6}$ for \nth{4}-order Runge Kutta. Red star denotes optimum, dashed line indicates optimal $t_0$ given $\dt$.}
    \label{fig:totalcost_dt_t0_emod_RK4}
\end{figure}
In Figure~\ref{fig:totalcost_dt_t0_emod_FE}, we show $\emodel$ for Forward Euler at a fixed cost of $U= \num{1.2e6}$ (the optimum is denoted by a red star). Moving to the right, discretization error becomes the dominant factor as $\dt \gg T_d$. The diagonal boundary gives the region of feasibility at which, under the cost constraint, sampling no longer occurs ($T_s= 0$). Moving from the optimum towards the bottom left, $t_0 \to 0$, $T_s \to 0$, and $\dt \ll T_d$; thus the transient error and sampling error become dominant. Similar plots for RK3 and RK4 are found in Figures~\ref{fig:totalcost_dt_t0_emod_RK3} and~\ref{fig:totalcost_dt_t0_emod_RK4}. The optimal errors and optimizing simulations are described in Table~\ref{tab:bigU_opt}. We can see from these results that, at a fixed budget with $U= \num{1.2e6}$, the effect of increasing the discretization order is make a smaller error possible with a larger timestep, which means fewer timesteps to traverse the spin-up time. These two effects combine to allow for an increase in the sampling time available $T_s$, allowing significantly less sampling error for RK3 compared to FE, and an additional-- albeit smaller-- benefit moving from RK3 to RK4, holding cost fixed.
\begin{table}[h]
    \centering
    \begin{tabular}{c|c|c|c|c|c}
        method &    $p$ &   $\emodel$ &     $\dt$ &         $t_0$ &     $T_s$ \\
        \hline
        FE &        $1$ &   \num{0.0502} &  \num{2.54e-4} & 30.2  &     275 \\
        RK3 &       $3$ &   \num{0.0130} &  \num{8.52e-3} & 35.5  &     3370 \\
        RK4 &       $4$ &   \num{8.89e-3} & \num{0.0224} &  36.7  &     6670
    \end{tabular}
    \caption{Optimal Lorenz simulations for output $g= u_2$ under budget of $U= \num{1.2e6}$ right-hand side evaluations using $\Mens= 1$.}
    \label{tab:bigU_opt}
\end{table}

\begin{figure}[h]
    \centering
    \includegraphics[width= \imgwidth]{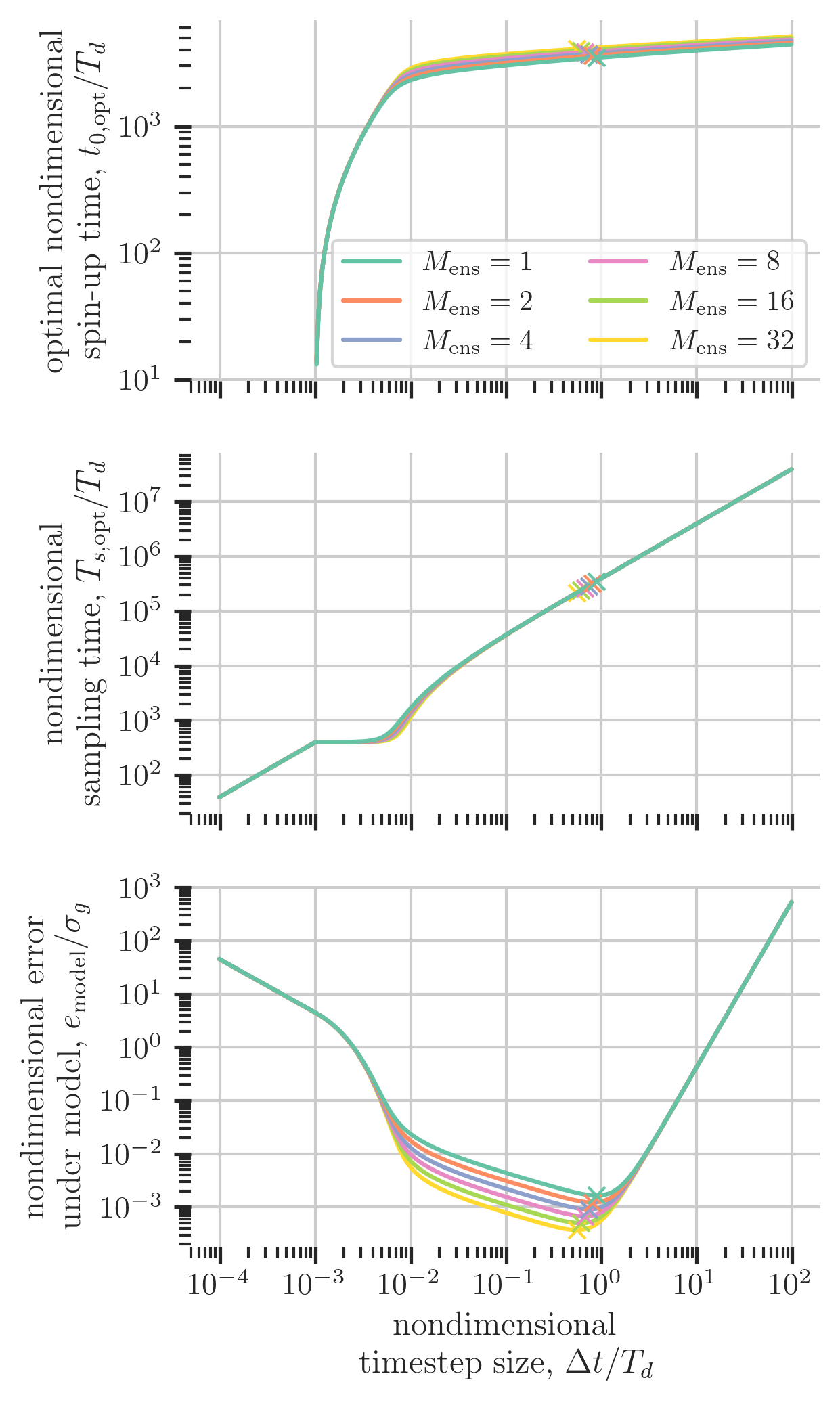}
    \caption{Dependence of normalized spin-up time $t_0/T_d$, sampling time $T_s/T_d$, and model error $\emodel/\sigma_g$ on normalized timestep $\dt/T_d$ with total cost set at $U= \num{1.2e6}$ for \nth{3}-order Runge Kutta.}
    \label{fig:totalcost_dt_combo_RK3}
\end{figure}

In Figure~\ref{fig:totalcost_dt_combo_RK3}, we take another perspective on these results for RK3 by varying $\dt$ and plotting the optimal $t_0$, $T_s$, and $\emodel$.
As $\dt$ gets large, the optimal choice of $t_0$ has logarithmic growth, and when $\dt/T_d \ll 1$, the optimal choice of $t_0$ rapidly falls to zero. Parallelization has a small but non-zero effect on the optimal choice of sample time.
The sampling time also has a small effect from parallelization, in this case constrained to a small region. Outside that $\dt$ region, $T_s$ scales with $\dt$ both as $\dt \to 0$ and as $\dt \to \infty$.

The bottom plot of Figure~\ref{fig:totalcost_dt_combo_RK3} shows the variation of error with $\dt$. In this plot we can see three distinct regions. For $\dt/T_d \gg 10^{-2}$, discretization error is the dominating error, and the convergence goes with the discretization error rate. Approaching the optimum, sampling error becomes the dominant error contribution, starting at $\dt \approx \num{2e-2}$ until $\dt \approx \num{e-3}$. In this region, the convergence is around the CLT-implied $1/2$ rate, and the effect of parallelization is clearly seen. For $\dt \lesssim \num{e-3}$, however, the spin-up error becomes the dominant error contribution. The optimal choice of $t_0$ begins to fall rapidly, as the sampling and spin-up must compete for computational resources under the budget. Once the spin-up error dominates, the paradigm by which \eqref{eq:abs_transient_model} is controlled shifts from the $\exp(-t_0)$ term to the $T_s^{-1}$ term as $\dt/T_d \to 0$, since resolving $T_s$ delivers both spin-up and sampling error control.

\begin{figure}[h]
    \centering
    \includegraphics[width= \imgwidth]{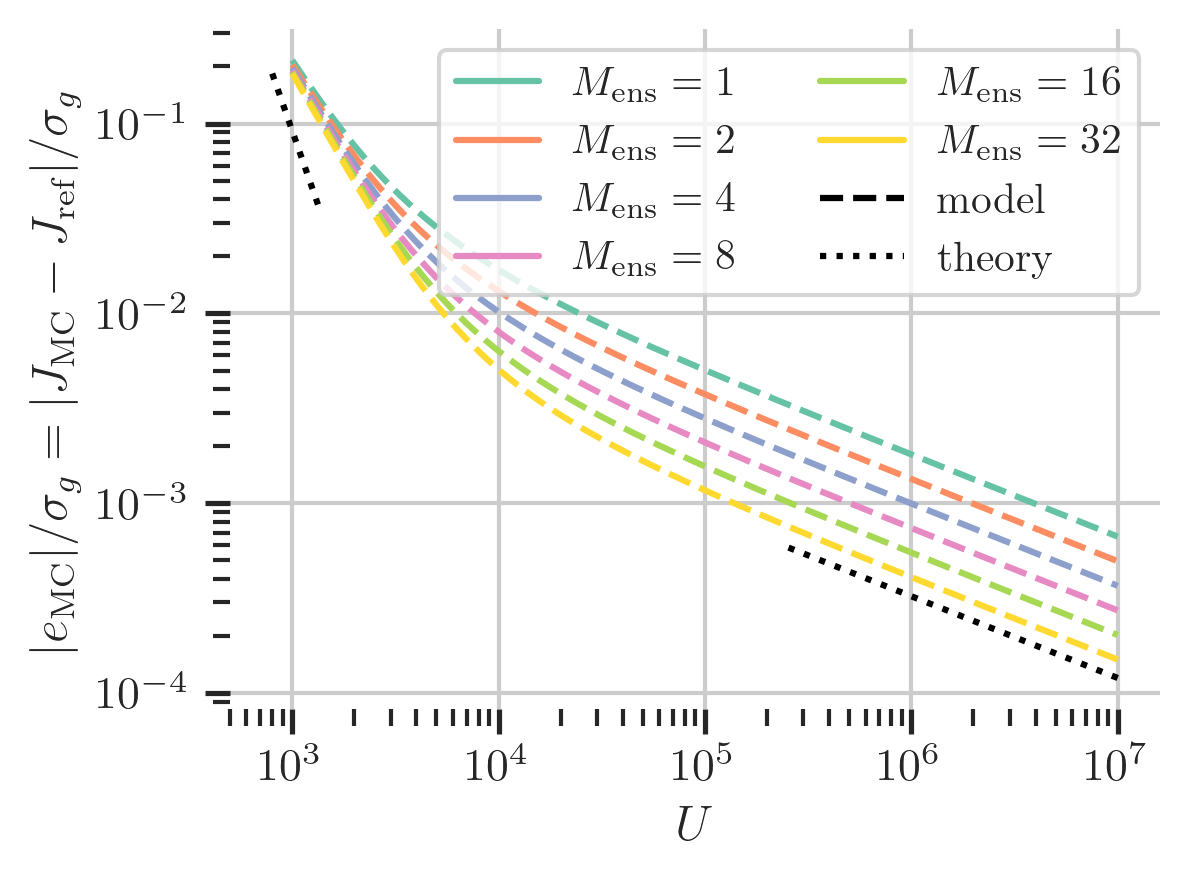}
    \caption{Optimal non-dimensional error under model as a function of total cost $U$ for RK3. Theory totem on left-hand side: discrete convergence rate, $1/q$ ; on right-hand side: $\frac{2(q+r)}{q}$ rate from \eqref{eq:error_optimal_mc_nondim}.}
    \label{fig:totalcost_scaling_model_RK3}
\end{figure}

This interdependence will evidently have an effect on the overall scaling between cost and error, which we now seek to understand. Here, we study the variation of $\eMC$ with $U$ under the optimal choices and evaluate how well $\eMC$ approximates experimental data for $\expect[|J_\mathrm{MC} - \Jinf|]$. In Figure~\ref{fig:totalcost_scaling_model_RK3}, the variation of $\eMC$ computed via \eqref{eq:error_model_constrained} as a function of $\Mens$ and $U$ is shown.
From this figure, we can see that, in the limit of small error, the sampling costs dominate and the best possible rate is given by the estimate in \eqref{eq:error_optimal_mc_nondim}, limited by the CLT. On the other hand, when the cost is more moderate, scaling of the error is close to the discretization error convergence rate in \eqref{eq:truncation_error}. In this region, the spin-up costs are significant, and high-order discretization brings the state more efficiently to the start of sampling. In the spin-up dominated region, the effect of the parallel ensemble approach is minimal since spin-up must be overcome on each processor.

\begin{figure}[h]
    \centering
    \includegraphics[width= \imgwidth]{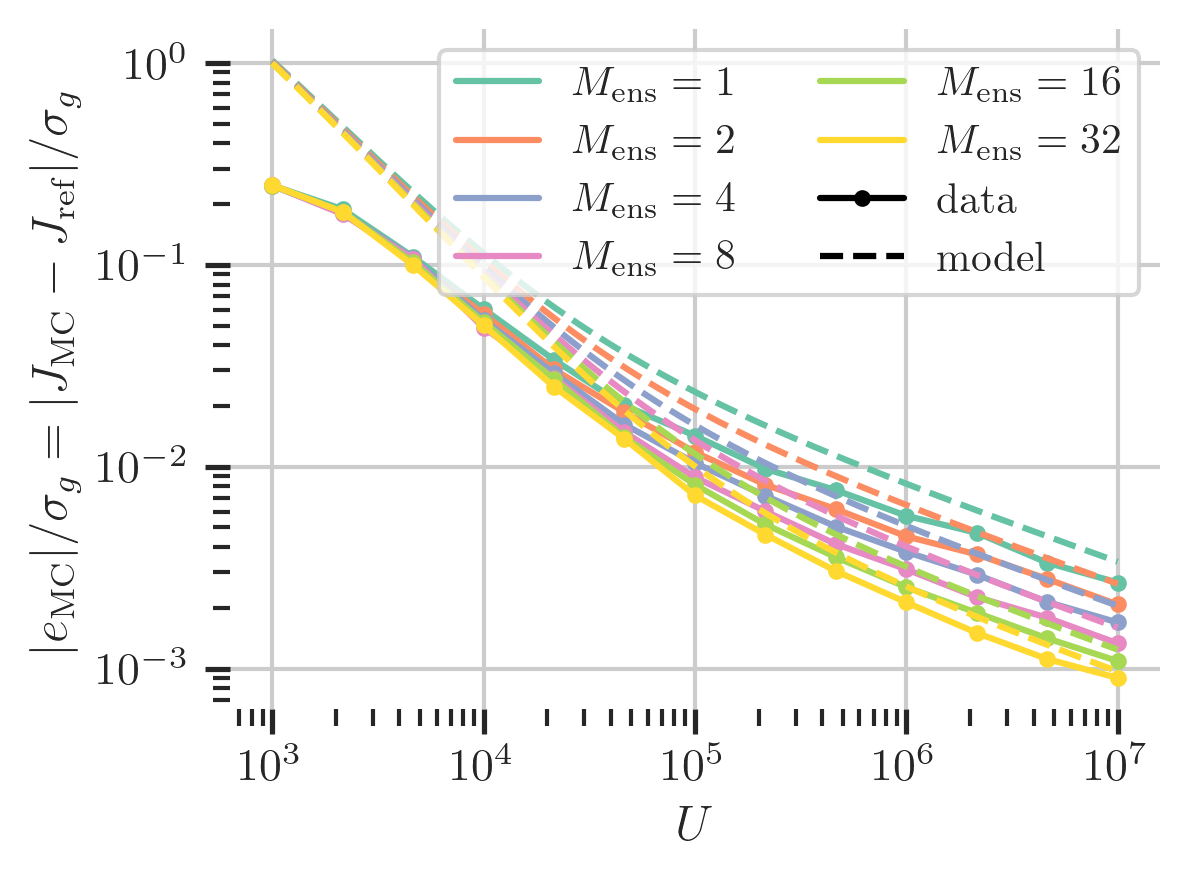}
    \caption{Total cost model and Monte Carlo validation as a function of total cost $U$ for FE.}
    \label{fig:totalcost_scaling_FE}
\end{figure}

\begin{figure}[h]
    \centering
    \includegraphics[width= \imgwidth]{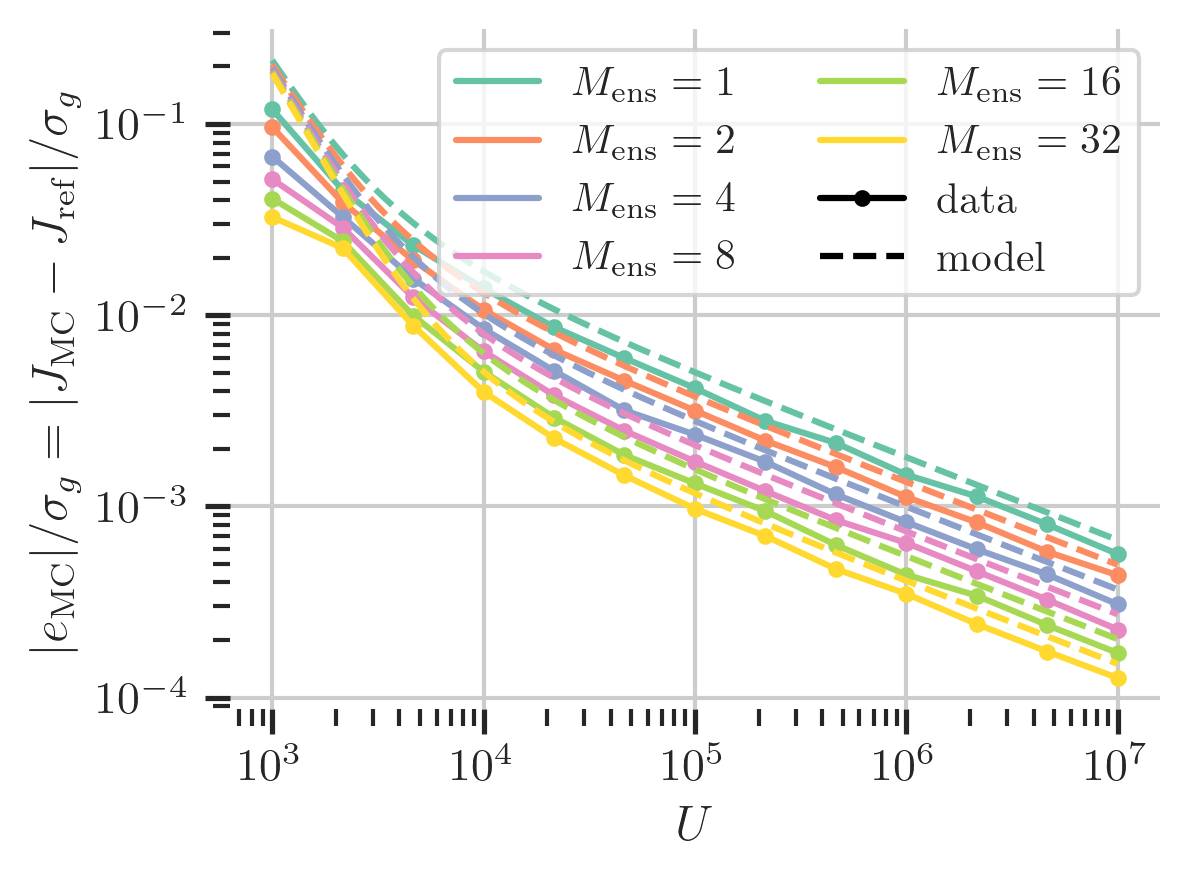}
    \caption{Total cost model and Monte Carlo validation as a function of total cost $U$ for RK3.}
    \label{fig:totalcost_scaling_RK3}
\end{figure}

\begin{figure}[h]
    \centering
    \includegraphics[width= \imgwidth]{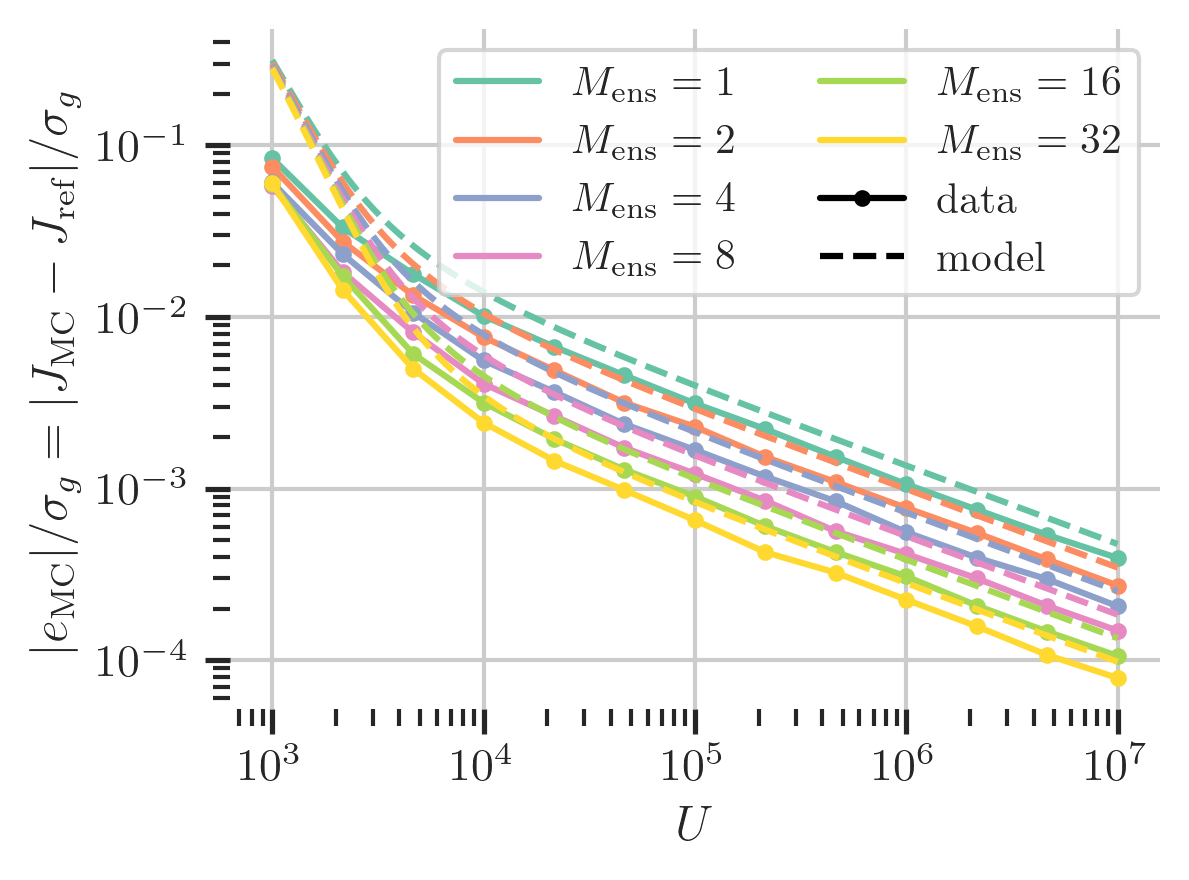}
    \caption{Total cost model and Monte Carlo validation as a function of total cost $U$ for RK4.}
    \label{fig:totalcost_scaling_RK4}
\end{figure}

Now, we validate the total error model for the Lorenz system by a final numerical experiment. At each choice of $\Mens$ and $U$, we generate \num{1000} individual realizations of $J_\mathrm{MC}$ at the computed $\dtMCopt$ and $\NMCopt$ and using the model fit given in Table~\ref{tab:NLS_fits_big}. In Figures~\ref{fig:totalcost_scaling_FE}, \ref{fig:totalcost_scaling_RK3}, and~\ref{fig:totalcost_scaling_RK4}, we show the predictions and the results of Monte Carlo estimates of $\expect[|J_\mathrm{MC} - \Jinf|]$ for our three discretizations.
These results validate the model, with significant discrepancies only when the asymptotic assumptions-- $\dt$ small and $T_s$ large-- do not hold, due to budget limitations in the limit of small $U$.


\section{Conclusions \& forthcoming work}

In this manuscript, we have developed a theoretical framework for the total error incurred by the discrete sampling of mean outputs of ergodic ODEs. These findings are validated by Monte Carlo studies of the Lorenz system using Runge-Kutta methods. We incorporate effects of parallelization and spin-up and validate that the models match observed results in experiments. Using these models, we are able to develop a comprehensive understanding of the relationship between the wall-clock cost of a simulation and the amount of error in expectation that it might achieve.

A key problem with the applicability of this research presented in this paper is the expense of identifying the parameters of the error model. In order to overcome this, we believe that leveraging a Bayesian approach as in \citet{oliver2014estimating} can allow us to approximate the model in \eqref{eq:error_model} at relatively small cost, and then exploit the result to conduct a high-fidelity simulation at (approximately) optimal discretizations. Further, the framework must be extended to handle chaotic PDE systems as opposed to ODE systems. Though many discrete PDE systems are discretized in a form that reduces to an ODE system, a rigorous model for the error and cost of a PDE system should account for the contributions of both temporal discretization \emph{and} spatial discretization. These will be the primary concerns of our forthcoming work.


\begin{acknowledgments}
The authors would like to acknowledge the support of The Boeing Company (technical monitor Dr. Andrew Cary).

\end{acknowledgments}

\section*{Conflicts of Interest}

The authors have no conflicts of interest to report.

\section*{Data Availability}

The data that support the findings of this study are available from the corresponding author upon reasonable request.


\bibliography{paper}

\begin{thebibliography}{29}%
\makeatletter
\providecommand \@ifxundefined [1]{%
 \@ifx{#1\undefined}
}%
\providecommand \@ifnum [1]{%
 \ifnum #1\expandafter \@firstoftwo
 \else \expandafter \@secondoftwo
 \fi
}%
\providecommand \@ifx [1]{%
 \ifx #1\expandafter \@firstoftwo
 \else \expandafter \@secondoftwo
 \fi
}%
\providecommand \natexlab [1]{#1}%
\providecommand \enquote  [1]{``#1''}%
\providecommand \bibnamefont  [1]{#1}%
\providecommand \bibfnamefont [1]{#1}%
\providecommand \citenamefont [1]{#1}%
\providecommand \href@noop [0]{\@secondoftwo}%
\providecommand \href [0]{\begingroup \@sanitize@url \@href}%
\providecommand \@href[1]{\@@startlink{#1}\@@href}%
\providecommand \@@href[1]{\endgroup#1\@@endlink}%
\providecommand \@sanitize@url [0]{\catcode `\\12\catcode `\$12\catcode
  `\&12\catcode `\#12\catcode `\^12\catcode `\_12\catcode `\%12\relax}%
\providecommand \@@startlink[1]{}%
\providecommand \@@endlink[0]{}%
\providecommand \url  [0]{\begingroup\@sanitize@url \@url }%
\providecommand \@url [1]{\endgroup\@href {#1}{\urlprefix }}%
\providecommand \urlprefix  [0]{URL }%
\providecommand \Eprint [0]{\href }%
\providecommand \doibase [0]{http://dx.doi.org/}%
\providecommand \selectlanguage [0]{\@gobble}%
\providecommand \bibinfo  [0]{\@secondoftwo}%
\providecommand \bibfield  [0]{\@secondoftwo}%
\providecommand \translation [1]{[#1]}%
\providecommand \BibitemOpen [0]{}%
\providecommand \bibitemStop [0]{}%
\providecommand \bibitemNoStop [0]{.\EOS\space}%
\providecommand \EOS [0]{\spacefactor3000\relax}%
\providecommand \BibitemShut  [1]{\csname bibitem#1\endcsname}%
\let\auto@bib@innerbib\@empty
\bibitem [{\citenamefont {Lighthill}(1986)}]{lighthill1986recently}%
  \BibitemOpen
  \bibfield  {author} {\bibinfo {author} {\bibfnamefont {M.~J.}\ \bibnamefont
  {Lighthill}},\ }\bibfield  {title} {\enquote {\bibinfo {title} {The recently
  recognized failure of predictability in {Newtonian} dynamics},}\ }\href@noop
  {} {\bibfield  {journal} {\bibinfo  {journal} {{Proceedings of the Royal
  Society of London. A. Mathematical and Physical Sciences}}\ }\textbf
  {\bibinfo {volume} {407}},\ \bibinfo {pages} {35--50} (\bibinfo {year}
  {1986})}\BibitemShut {NoStop}%
\bibitem [{\citenamefont {Eckmann}\ and\ \citenamefont
  {Ruelle}(1985)}]{eckmann1985ergodic}%
  \BibitemOpen
  \bibfield  {author} {\bibinfo {author} {\bibfnamefont {J.-P.}\ \bibnamefont
  {Eckmann}}\ and\ \bibinfo {author} {\bibfnamefont {D.}~\bibnamefont
  {Ruelle}},\ }\bibfield  {title} {\enquote {\bibinfo {title} {Ergodic theory
  of chaos and strange attractors},}\ }in\ \href@noop {} {\emph {\bibinfo
  {booktitle} {The theory of chaotic attractors}}}\ (\bibinfo  {publisher}
  {Springer},\ \bibinfo {year} {1985})\ pp.\ \bibinfo {pages}
  {273--312}\BibitemShut {NoStop}%
\bibitem [{\citenamefont {Chapman}(1979)}]{chapman1979computational}%
  \BibitemOpen
  \bibfield  {author} {\bibinfo {author} {\bibfnamefont {D.~R.}\ \bibnamefont
  {Chapman}},\ }\bibfield  {title} {\enquote {\bibinfo {title} {Computational
  aerodynamics development and outlook},}\ }\href@noop {} {\bibfield  {journal}
  {\bibinfo  {journal} {AIAA journal}\ }\textbf {\bibinfo {volume} {17}},\
  \bibinfo {pages} {1293--1313} (\bibinfo {year} {1979})}\BibitemShut {NoStop}%
\bibitem [{\citenamefont {Spalart}\ \emph {et~al.}(1997)\citenamefont
  {Spalart}, \citenamefont {Jou}, \citenamefont {Strelets}, \citenamefont
  {Allmaras} \emph {et~al.}}]{spalart1997comments}%
  \BibitemOpen
  \bibfield  {author} {\bibinfo {author} {\bibfnamefont {P.}~\bibnamefont
  {Spalart}}, \bibinfo {author} {\bibfnamefont {W.}~\bibnamefont {Jou}},
  \bibinfo {author} {\bibfnamefont {M.}~\bibnamefont {Strelets}}, \bibinfo
  {author} {\bibfnamefont {S.}~\bibnamefont {Allmaras}},  \emph {et~al.},\
  }\bibfield  {title} {\enquote {\bibinfo {title} {{Comments on the feasibility
  of LES for wings, and on a hybrid RANS/LES approach}},}\ }\href@noop {}
  {\bibfield  {journal} {\bibinfo  {journal} {Advances in DNS/LES}\ }\textbf
  {\bibinfo {volume} {1}},\ \bibinfo {pages} {4--8} (\bibinfo {year}
  {1997})}\BibitemShut {NoStop}%
\bibitem [{\citenamefont {Choi}\ and\ \citenamefont
  {Moin}(2012)}]{choi2012grid}%
  \BibitemOpen
  \bibfield  {author} {\bibinfo {author} {\bibfnamefont {H.}~\bibnamefont
  {Choi}}\ and\ \bibinfo {author} {\bibfnamefont {P.}~\bibnamefont {Moin}},\
  }\bibfield  {title} {\enquote {\bibinfo {title} {Grid-point requirements for
  large eddy simulation: {Chapman}’s estimates revisited},}\ }\href@noop {}
  {\bibfield  {journal} {\bibinfo  {journal} {Physics of Fluids}\ }\textbf
  {\bibinfo {volume} {24}},\ \bibinfo {pages} {011702} (\bibinfo {year}
  {2012})}\BibitemShut {NoStop}%
\bibitem [{\citenamefont {Kim}, \citenamefont {Moin},\ and\ \citenamefont
  {Moser}(1987)}]{kim1987turbulence}%
  \BibitemOpen
  \bibfield  {author} {\bibinfo {author} {\bibfnamefont {J.}~\bibnamefont
  {Kim}}, \bibinfo {author} {\bibfnamefont {P.}~\bibnamefont {Moin}}, \ and\
  \bibinfo {author} {\bibfnamefont {R.}~\bibnamefont {Moser}},\ }\bibfield
  {title} {\enquote {\bibinfo {title} {Turbulence statistics in fully developed
  channel flow at low {Reynolds} number},}\ }\href {\doibase
  10.1017/S0022112087000892} {\bibfield  {journal} {\bibinfo  {journal}
  {Journal of Fluid Mechanics}\ }\textbf {\bibinfo {volume} {177}},\ \bibinfo
  {pages} {133–166} (\bibinfo {year} {1987})}\BibitemShut {NoStop}%
\bibitem [{\citenamefont {Lozano-Dur\'{a}n}\ and\ \citenamefont
  {Jim\'{e}nez}(2014)}]{lozano2014effect}%
  \BibitemOpen
  \bibfield  {author} {\bibinfo {author} {\bibfnamefont {A.}~\bibnamefont
  {Lozano-Dur\'{a}n}}\ and\ \bibinfo {author} {\bibfnamefont {J.}~\bibnamefont
  {Jim\'{e}nez}},\ }\bibfield  {title} {\enquote {\bibinfo {title} {Effect of
  the computational domain on direct simulations of turbulent channels up to
  $\mathrm{Re}_{\tau}= 4200$},}\ }\href {\doibase 10.1063/1.4862918} {\bibfield
   {journal} {\bibinfo  {journal} {Physics of Fluids}\ }\textbf {\bibinfo
  {volume} {26}},\ \bibinfo {pages} {011702} (\bibinfo {year} {2014})},\
  \Eprint {http://arxiv.org/abs/https://doi.org/10.1063/1.4862918}
  {https://doi.org/10.1063/1.4862918} \BibitemShut {NoStop}%
\bibitem [{\citenamefont {Del~Álamo}\ \emph {et~al.}(2004)\citenamefont
  {Del~Álamo}, \citenamefont {Jim\'{e}nez}, \citenamefont {Zandonade},\ and\
  \citenamefont {Moser}}]{delalamo2004scaling}%
  \BibitemOpen
  \bibfield  {author} {\bibinfo {author} {\bibfnamefont {J.~C.}\ \bibnamefont
  {Del~Álamo}}, \bibinfo {author} {\bibfnamefont {J.}~\bibnamefont
  {Jim\'{e}nez}}, \bibinfo {author} {\bibfnamefont {P.}~\bibnamefont
  {Zandonade}}, \ and\ \bibinfo {author} {\bibfnamefont {R.~D.}\ \bibnamefont
  {Moser}},\ }\bibfield  {title} {\enquote {\bibinfo {title} {Scaling of the
  energy spectra of turbulent channels},}\ }\href {\doibase
  10.1017/S002211200300733X} {\bibfield  {journal} {\bibinfo  {journal}
  {Journal of Fluid Mechanics}\ }\textbf {\bibinfo {volume} {500}},\ \bibinfo
  {pages} {135–144} (\bibinfo {year} {2004})}\BibitemShut {NoStop}%
\bibitem [{\citenamefont {Goc}\ \emph {et~al.}(2021)\citenamefont {Goc},
  \citenamefont {Lehmkuhl}, \citenamefont {Park}, \citenamefont {Bose},\ and\
  \citenamefont {Moin}}]{goc2021large}%
  \BibitemOpen
  \bibfield  {author} {\bibinfo {author} {\bibfnamefont {K.~A.}\ \bibnamefont
  {Goc}}, \bibinfo {author} {\bibfnamefont {O.}~\bibnamefont {Lehmkuhl}},
  \bibinfo {author} {\bibfnamefont {G.~I.}\ \bibnamefont {Park}}, \bibinfo
  {author} {\bibfnamefont {S.~T.}\ \bibnamefont {Bose}}, \ and\ \bibinfo
  {author} {\bibfnamefont {P.}~\bibnamefont {Moin}},\ }\bibfield  {title}
  {\enquote {\bibinfo {title} {Large eddy simulation of aircraft at affordable
  cost: a milestone in computational fluid dynamics},}\ }\href {\doibase
  10.1017/flo.2021.17} {\bibfield  {journal} {\bibinfo  {journal} {Flow}\
  }\textbf {\bibinfo {volume} {1}},\ \bibinfo {pages} {E14} (\bibinfo {year}
  {2021})}\BibitemShut {NoStop}%
\bibitem [{\citenamefont {Thompson}\ \emph {et~al.}(2016)\citenamefont
  {Thompson}, \citenamefont {Sampaio}, \citenamefont {de~Bragan{\c{c}}a~Alves},
  \citenamefont {Thais},\ and\ \citenamefont
  {Mompean}}]{thompson2016methodology}%
  \BibitemOpen
  \bibfield  {author} {\bibinfo {author} {\bibfnamefont {R.~L.}\ \bibnamefont
  {Thompson}}, \bibinfo {author} {\bibfnamefont {L.~E.~B.}\ \bibnamefont
  {Sampaio}}, \bibinfo {author} {\bibfnamefont {F.~A.}\ \bibnamefont
  {de~Bragan{\c{c}}a~Alves}}, \bibinfo {author} {\bibfnamefont
  {L.}~\bibnamefont {Thais}}, \ and\ \bibinfo {author} {\bibfnamefont
  {G.}~\bibnamefont {Mompean}},\ }\bibfield  {title} {\enquote {\bibinfo
  {title} {A methodology to evaluate statistical errors in {DNS} data of plane
  channel flows},}\ }\href@noop {} {\bibfield  {journal} {\bibinfo  {journal}
  {Computers \& Fluids}\ }\textbf {\bibinfo {volume} {130}},\ \bibinfo {pages}
  {1--7} (\bibinfo {year} {2016})}\BibitemShut {NoStop}%
\bibitem [{\citenamefont {Russo}\ and\ \citenamefont
  {Luchini}(2017)}]{russo2017fast}%
  \BibitemOpen
  \bibfield  {author} {\bibinfo {author} {\bibfnamefont {S.}~\bibnamefont
  {Russo}}\ and\ \bibinfo {author} {\bibfnamefont {P.}~\bibnamefont
  {Luchini}},\ }\bibfield  {title} {\enquote {\bibinfo {title} {A fast
  algorithm for the estimation of statistical error in {DNS} (or experimental)
  time averages},}\ }\href@noop {} {\bibfield  {journal} {\bibinfo  {journal}
  {Journal of Computational Physics}\ }\textbf {\bibinfo {volume} {347}},\
  \bibinfo {pages} {328--340} (\bibinfo {year} {2017})}\BibitemShut {NoStop}%
\bibitem [{\citenamefont {Mockett}, \citenamefont {Knacke},\ and\ \citenamefont
  {Thiele}(2010)}]{mockett2010detection}%
  \BibitemOpen
  \bibfield  {author} {\bibinfo {author} {\bibfnamefont {C.}~\bibnamefont
  {Mockett}}, \bibinfo {author} {\bibfnamefont {T.}~\bibnamefont {Knacke}}, \
  and\ \bibinfo {author} {\bibfnamefont {F.}~\bibnamefont {Thiele}},\
  }\bibfield  {title} {\enquote {\bibinfo {title} {Detection of initial
  transient and estimation of statistical error in time-resolved turbulent flow
  data},}\ }in\ \href@noop {} {\emph {\bibinfo {booktitle} {Proceedings of the
  8th International Symposium on Engineering Turbulence Modelling and
  Measurements}}}\ (\bibinfo {organization} {European Research Collaboration on
  Flow Turbulence and Combustion},\ \bibinfo {year} {2010})\ pp.\ \bibinfo
  {pages} {9--11}\BibitemShut {NoStop}%
\bibitem [{\citenamefont {Oliver}\ \emph {et~al.}(2014)\citenamefont {Oliver},
  \citenamefont {Malaya}, \citenamefont {Ulerich},\ and\ \citenamefont
  {Moser}}]{oliver2014estimating}%
  \BibitemOpen
  \bibfield  {author} {\bibinfo {author} {\bibfnamefont {T.~A.}\ \bibnamefont
  {Oliver}}, \bibinfo {author} {\bibfnamefont {N.}~\bibnamefont {Malaya}},
  \bibinfo {author} {\bibfnamefont {R.}~\bibnamefont {Ulerich}}, \ and\
  \bibinfo {author} {\bibfnamefont {R.~D.}\ \bibnamefont {Moser}},\ }\bibfield
  {title} {\enquote {\bibinfo {title} {Estimating uncertainties in statistics
  computed from direct numerical simulation},}\ }\href@noop {} {\bibfield
  {journal} {\bibinfo  {journal} {Physics of Fluids}\ }\textbf {\bibinfo
  {volume} {26}},\ \bibinfo {pages} {035101} (\bibinfo {year}
  {2014})}\BibitemShut {NoStop}%
\bibitem [{\citenamefont {Lee}\ and\ \citenamefont
  {Moser}(2015)}]{lee2015direct}%
  \BibitemOpen
  \bibfield  {author} {\bibinfo {author} {\bibfnamefont {M.}~\bibnamefont
  {Lee}}\ and\ \bibinfo {author} {\bibfnamefont {R.~D.}\ \bibnamefont
  {Moser}},\ }\bibfield  {title} {\enquote {\bibinfo {title} {Direct numerical
  simulation of turbulent channel flow up to $\mathrm{Re}_{\tau} \approx
  5200$},}\ }\href {\doibase 10.1017/jfm.2015.268} {\bibfield  {journal}
  {\bibinfo  {journal} {Journal of Fluid Mechanics}\ }\textbf {\bibinfo
  {volume} {774}},\ \bibinfo {pages} {395–415} (\bibinfo {year}
  {2015})}\BibitemShut {NoStop}%
\bibitem [{\citenamefont {Stuart}(1994)}]{stuart1994numerical}%
  \BibitemOpen
  \bibfield  {author} {\bibinfo {author} {\bibfnamefont {A.~M.}\ \bibnamefont
  {Stuart}},\ }\bibfield  {title} {\enquote {\bibinfo {title} {Numerical
  analysis of dynamical systems},}\ }\href@noop {} {\bibfield  {journal}
  {\bibinfo  {journal} {Acta numerica}\ }\textbf {\bibinfo {volume} {3}},\
  \bibinfo {pages} {467--572} (\bibinfo {year} {1994})}\BibitemShut {NoStop}%
\bibitem [{\citenamefont {Denker}(1989)}]{denker1989central}%
  \BibitemOpen
  \bibfield  {author} {\bibinfo {author} {\bibfnamefont {M.}~\bibnamefont
  {Denker}},\ }\bibfield  {title} {\enquote {\bibinfo {title} {The central
  limit theorem for dynamical systems},}\ }\href@noop {} {\bibfield  {journal}
  {\bibinfo  {journal} {Banach Center Publications}\ }\textbf {\bibinfo
  {volume} {1}},\ \bibinfo {pages} {33--62} (\bibinfo {year}
  {1989})}\BibitemShut {NoStop}%
\bibitem [{\citenamefont {Bradley}(2005)}]{bradley2005basic}%
  \BibitemOpen
  \bibfield  {author} {\bibinfo {author} {\bibfnamefont {R.~C.}\ \bibnamefont
  {Bradley}},\ }\bibfield  {title} {\enquote {\bibinfo {title} {Basic
  properties of strong mixing conditions. a survey and some open questions},}\
  }\href@noop {} {\bibfield  {journal} {\bibinfo  {journal} {Probability
  Surveys}\ }\textbf {\bibinfo {volume} {2}},\ \bibinfo {pages} {107--144}
  (\bibinfo {year} {2005})}\BibitemShut {NoStop}%
\bibitem [{\citenamefont {Ara{\'u}jo}, \citenamefont {Melbourne},\ and\
  \citenamefont {Varandas}(2015)}]{araujo2015rapid}%
  \BibitemOpen
  \bibfield  {author} {\bibinfo {author} {\bibfnamefont {V.}~\bibnamefont
  {Ara{\'u}jo}}, \bibinfo {author} {\bibfnamefont {I.}~\bibnamefont
  {Melbourne}}, \ and\ \bibinfo {author} {\bibfnamefont {P.}~\bibnamefont
  {Varandas}},\ }\bibfield  {title} {\enquote {\bibinfo {title} {Rapid mixing
  for the {Lorenz} attractor and statistical limit laws for their time-1
  maps},}\ }\href@noop {} {\bibfield  {journal} {\bibinfo  {journal}
  {Communications in Mathematical Physics}\ }\textbf {\bibinfo {volume}
  {340}},\ \bibinfo {pages} {901--938} (\bibinfo {year} {2015})}\BibitemShut
  {NoStop}%
\bibitem [{\citenamefont {Hairer}(1993)}]{hairer1993solving}%
  \BibitemOpen
  \bibfield  {author} {\bibinfo {author} {\bibfnamefont {E.}~\bibnamefont
  {Hairer}},\ }\href@noop {} {\emph {\bibinfo {title} {Solving ordinary
  differential equations {II}: stiff and differential-algebraic problems}}},\
  \bibinfo {edition} {second edition}\ ed.,\ \bibinfo {series} {{Springer
  Series in Computational Mathematics}}\ No.~\bibinfo {number} {14}\ (\bibinfo
  {publisher} {Springer},\ \bibinfo {address} {Berlin, Germany},\ \bibinfo
  {year} {1993})\BibitemShut {NoStop}%
\bibitem [{\citenamefont {Viswanath}(2001)}]{viswanath2001global}%
  \BibitemOpen
  \bibfield  {author} {\bibinfo {author} {\bibfnamefont {D.}~\bibnamefont
  {Viswanath}},\ }\bibfield  {title} {\enquote {\bibinfo {title} {Global errors
  of numerical {ODE} solvers and {Lyapunov's} theory of stability},}\
  }\href@noop {} {\bibfield  {journal} {\bibinfo  {journal} {{IMA} Journal of
  Numerical Analysis}\ }\textbf {\bibinfo {volume} {21}},\ \bibinfo {pages}
  {387--406} (\bibinfo {year} {2001})}\BibitemShut {NoStop}%
\bibitem [{\citenamefont {Lorenz}(1963)}]{lorenz1963deterministic}%
  \BibitemOpen
  \bibfield  {author} {\bibinfo {author} {\bibfnamefont {E.~N.}\ \bibnamefont
  {Lorenz}},\ }\bibfield  {title} {\enquote {\bibinfo {title} {Deterministic
  nonperiodic flow},}\ }\href@noop {} {\bibfield  {journal} {\bibinfo
  {journal} {Journal of the Atmospheric Sciences}\ }\textbf {\bibinfo {volume}
  {20}},\ \bibinfo {pages} {130--141} (\bibinfo {year} {1963})}\BibitemShut
  {NoStop}%
\bibitem [{\citenamefont {Sparrow}(1982)}]{sparrow1982lorenz}%
  \BibitemOpen
  \bibfield  {author} {\bibinfo {author} {\bibfnamefont {C.}~\bibnamefont
  {Sparrow}},\ }\href@noop {} {\emph {\bibinfo {title} {The Lorenz equations:
  bifurcations, chaos, and strange attractors}}}\ (\bibinfo  {publisher}
  {Springer Science \& Business Media},\ \bibinfo {year} {1982})\BibitemShut
  {NoStop}%
\bibitem [{\citenamefont {Dormand}, \citenamefont {Duckers},\ and\
  \citenamefont {Prince}(1984)}]{dormand1984global}%
  \BibitemOpen
  \bibfield  {author} {\bibinfo {author} {\bibfnamefont {J.}~\bibnamefont
  {Dormand}}, \bibinfo {author} {\bibfnamefont {R.}~\bibnamefont {Duckers}}, \
  and\ \bibinfo {author} {\bibfnamefont {P.}~\bibnamefont {Prince}},\
  }\bibfield  {title} {\enquote {\bibinfo {title} {Global error estimation with
  {Runge-Kutta} methods},}\ }\href@noop {} {\bibfield  {journal} {\bibinfo
  {journal} {IMA Journal of Numerical Analysis}\ }\textbf {\bibinfo {volume}
  {4}},\ \bibinfo {pages} {169--184} (\bibinfo {year} {1984})}\BibitemShut
  {NoStop}%
\bibitem [{\citenamefont {Trenberth}(1984)}]{trenberth1984some}%
  \BibitemOpen
  \bibfield  {author} {\bibinfo {author} {\bibfnamefont {K.~E.}\ \bibnamefont
  {Trenberth}},\ }\bibfield  {title} {\enquote {\bibinfo {title} {Some effects
  of finite sample size and persistence on meteorological statistics. {Part I}:
  {Autocorrelations}},}\ }\href@noop {} {\bibfield  {journal} {\bibinfo
  {journal} {Monthly Weather Review}\ }\textbf {\bibinfo {volume} {112}},\
  \bibinfo {pages} {2359--2368} (\bibinfo {year} {1984})}\BibitemShut {NoStop}%
\bibitem [{Note1()}]{Note1}%
  \BibitemOpen
  \bibinfo {note} {In general, we expect $q= p$, but due to cancellation of
  local errors, $q > p$ occurs in practice for the Lorenz system. In the
  expected case of $q= p$, we should expect $\protect \mathcal {G}_{\protect
  \mathrm {max}}= C_q/(c_p T_d)$.}\BibitemShut {Stop}%
\bibitem [{Note2()}]{Note2}%
  \BibitemOpen
  \bibinfo {note} {Derivatives of $f$ are computed analytically using the chain
  rule.}\BibitemShut {Stop}%
\bibitem [{\citenamefont {Harris}(1978)}]{harris1978windows}%
  \BibitemOpen
  \bibfield  {author} {\bibinfo {author} {\bibfnamefont {F.}~\bibnamefont
  {Harris}},\ }\bibfield  {title} {\enquote {\bibinfo {title} {On then use of
  windows for harmonic analysis with the discrete fourier transform},}\
  }\href@noop {} {\bibfield  {journal} {\bibinfo  {journal} {Proceedings of the
  IEEE}\ }\textbf {\bibinfo {volume} {60}} (\bibinfo {year}
  {1978})}\BibitemShut {NoStop}%
\bibitem [{\citenamefont {Karniadakis}\ and\ \citenamefont
  {Sherwin}(2005)}]{karniadakis2005spectral}%
  \BibitemOpen
  \bibfield  {author} {\bibinfo {author} {\bibfnamefont {G.}~\bibnamefont
  {Karniadakis}}\ and\ \bibinfo {author} {\bibfnamefont {S.}~\bibnamefont
  {Sherwin}},\ }\href@noop {} {\emph {\bibinfo {title} {Spectral/hp-element
  methods for computational fluid dynamics}}}\ (\bibinfo  {publisher} {Oxford
  University Press},\ \bibinfo {year} {2005})\BibitemShut {NoStop}%
\bibitem [{\citenamefont {Makarashvili}\ \emph {et~al.}(2017)\citenamefont
  {Makarashvili}, \citenamefont {Merzari}, \citenamefont {Obabko},
  \citenamefont {Siegel},\ and\ \citenamefont
  {Fischer}}]{makarashvili2017performance}%
  \BibitemOpen
  \bibfield  {author} {\bibinfo {author} {\bibfnamefont {V.}~\bibnamefont
  {Makarashvili}}, \bibinfo {author} {\bibfnamefont {E.}~\bibnamefont
  {Merzari}}, \bibinfo {author} {\bibfnamefont {A.}~\bibnamefont {Obabko}},
  \bibinfo {author} {\bibfnamefont {A.}~\bibnamefont {Siegel}}, \ and\ \bibinfo
  {author} {\bibfnamefont {P.}~\bibnamefont {Fischer}},\ }\bibfield  {title}
  {\enquote {\bibinfo {title} {A performance analysis of ensemble averaging for
  high fidelity turbulence simulations at the strong scaling limit},}\
  }\href@noop {} {\bibfield  {journal} {\bibinfo  {journal} {Computer Physics
  Communications}\ }\textbf {\bibinfo {volume} {219}},\ \bibinfo {pages}
  {236--245} (\bibinfo {year} {2017})}\BibitemShut {NoStop}%
\end{thebibliography}%

\end{document}